\documentclass[leqno, 12pt]{article}

\usepackage{makeidx}
\makeindex

\def\im{\mathop{\rm Im}}
\def\re{\mathop{\rm Re}}
\def\diam{\mathop{\rm diam}}

\newtheorem{theorem}{Theorem}
\newtheorem{lemma}[theorem]{Lemma}
\newtheorem{proposition}[theorem]{Proposition}
\newtheorem{definition}[theorem]{Definition}
\newtheorem{corollary}[theorem]{Corollary}

\newcommand{\begintheorem}{\addtocounter{equation}{1}\begin{theorem}}
\newcommand{\beginlemma}{\addtocounter{equation}{1}\begin{lemma}}
\newcommand{\beginproposition}{\addtocounter{equation}{1}\begin{proposition}}
\newcommand{\begindefinition}{\addtocounter{equation}{1}\begin{definition}}
\newcommand{\begincorollary}{\addtocounter{equation}{1}\begin{corollary}}

\renewcommand{\thetheorem}{\arabic{section}.\arabic{equation}}
\renewcommand{\theequation}{\arabic{section}.\arabic{equation}}

\begin{document}

\title{Adventures in Harmonic Analysis}

\author{Stephen William Semmes	\\
	Rice University		\\
	Houston, Texas}

\date{}

\maketitle

\tableofcontents

\bigskip

	As usual, the integers are denoted ${\bf Z}$,\index{#Z@${\bf
Z}$} and the real and complex numbers are denoted ${\bf
R}$,\index{#R@${\bf R}$} ${\bf C}$,\index{#C@${\bf C}$} respectively.
If $x$ is a real number, then its \emph{absolute value}\index{absolute
value} is denoted $|x|$ and defined to be equal to $x$ when $x \ge 0$
and to $-x$ when $x \le 0$.  Thus $|x| \ge 0$ for every $x \in {\bf
R}$ and $|x| = 0$ if and only if $x = 0$.  One can check that
\begin{equation}
	|x + y| \le |x| + |y|
\end{equation}
and
\begin{equation}
	|x \, y| = |x| \, |y|
\end{equation}
for all $x, y \in {\bf R}$.

	A complex number $z$ can be expressed in a unique way as $x +
y \, i$, where $x$, $y$ are real numbers and $i^2 = -1$.  We may refer
to $x$ and $y$ as the \emph{real} and \emph{imaginary parts} of $z$,
denoted $\re z$, $\im z$, respectively.  The \emph{complex
conjugate}\index{complex conjugation} of $z$ is denoted $\overline{z}$
and defined by
\begin{equation}
	\overline{z} = x - y \, i.
\end{equation}
Thus
\begin{equation}
	\re z = \frac{z + \overline{z}}{2}
\end{equation}
and
\begin{equation}
	\im z = \frac{z - \overline{z}}{2 \, i}.
\end{equation}
Also, for all $z, w \in {\bf C}$,
\begin{equation}
	\overline{z + w} = \overline{z} + \overline{w}
\end{equation}
and
\begin{equation}
	\overline{z \, w} = \overline{z} \, \overline{w}.
\end{equation}

	If $z = x + y \, i$ is a complex number, $x, y \in {\bf R}$,
then $z \, \overline{z} = x^2 + y^2$.  The
\emph{modulus}\index{modulus} of $z$ is denoted $|z|$ and is defined
by
\begin{equation}
	|z| = \sqrt{x^2 + y^2}.
\end{equation}
In particular, $|\re z|, |\im z| \le |z|$.  Because $|z \, w|^2 = z \,
\overline{z} \, w \, \overline{w}$, we get that
\begin{equation}
	|z \, w| = |z| \, |w|
\end{equation}
for every $z, w \in {\bf C}$.  Also,
\begin{equation}
	|z + w| \le |z| + |w|
\end{equation}
for every $z, w \in {\bf C}$, because
\begin{eqnarray}
	|z + w|^2 & = & (z + w) (\overline{z} + \overline{w})	\\
		& = & |z|^2 + 2 \re z \, \overline{w} + |w|^2	\nonumber \\
		& \le & |z|^2 + 2 \, |z| \, |w| + |w|^2		\nonumber \\
		& = & (|z| + |w|)^2.				\nonumber
\end{eqnarray}

\section{Power Series}
\label{power series}
\setcounter{equation}{0}

	By a \emph{power series}\index{power series} in one complex
variable we mean a series
\begin{equation}
\label{sum_{n = 0}^infty a_n z^n}
	\sum_{n = 0}^\infty a_n \, z^n
\end{equation}
where the coefficients $a_n$ are complex numbers, $z$ is a complex
variable, and $z^n$ is interpreted as being equal to $1$ for every $z
\in {\bf C}$ when $n = 0$.  Such a series converges trivially when $z
= 0$, and may or may not converge elsewhere.

	For example, consider the \emph{geometric
series}\index{geometric series}
\begin{equation}
	\sum_{n=0}^\infty z^n.
\end{equation}
For every $z \in {\bf C}$ and nonnegative integer $r$ we have that
\begin{equation}
	(1 - z) \sum_{n=0}^r z^n = 1 - z^{r+1},
\end{equation}
and hence
\begin{equation}
	\sum_{n=0}^r z^n = \frac{1 - z^{r+1}}{1 - z}
\end{equation}
when $z \ne 1$.  If $|z| < 1$, then the geometric series converges and
\begin{equation}
	\sum_{n=0}^\infty z^n = \frac{1}{1-z},
\end{equation}
while if $|z| \ge 1$ then $|z|^n \ge 1$ for all $n$ and the series
diverges because the terms do not tend to $0$ as $n \to \infty$.

	If (\ref{sum_{n = 0}^infty a_n z^n}) converges for some $z_0
\in {\bf C}$, then
\begin{equation}
	\lim_{n \to \infty} a_n \, z_0^n = 0.
\end{equation}
In particular, $\{a_n \, z_0^n\}_{n=0}^\infty$ is a bounded sequence
of complex numbers, which is to say that there is an $A \ge 0$ such
that
\begin{equation}
	|a_n| \, |z_0|^n \le A
\end{equation}
for every positive integer $n$.  This implies that (\ref{sum_{n =
0}^infty a_n z^n}) converges absolutely when $|z| < |z_0|$, by
comparison with a geometric series.  Moreover, the partial sums
converge uniformly on the set of $z \in {\bf C}$ with $|z| \le r$ for
every $r \ge 0$ such that $r < |z_0|$, as a consequence of the
Weierstrass $M$-test.  The \emph{radius of convergence}\index{radius
of convergence} of (\ref{sum_{n = 0}^infty a_n z^n}) is the unique
$R$, $0 \le R \le +\infty$, such that (\ref{sum_{n = 0}^infty a_n
z^n}) converges absolutely when $|z| < R$ and does not converge when
$|z| > R$.

	If (\ref{sum_{n = 0}^infty a_n z^n}) has radius of convergence
$R$, then the power series defines a continuous complex-valued
function on the open disk
\begin{equation}
	\{z \in {\bf C} : |z| < R \},
\end{equation}
which is the whole complex plane when $R = +\infty$.  If $R < +\infty$
and
\begin{equation}
	\sum_{n = 0}^\infty |a_n| \, R^n
\end{equation}
converges, then (\ref{sum_{n = 0}^infty a_n z^n}) converges absolutely
for every $z \in {\bf C}$ such that $|z| \le R$, the partial sums
of the series converges uniformly on the closed disk
\begin{equation}
	\{z \in {\bf C} : |z| \le R \},
\end{equation}
and (\ref{sum_{n = 0}^infty a_n z^n}) defines a continuous function on
this disk.

	One can multiply a pair of power series $\sum_{j = 0}^\infty
a_j \, z^j$, $\sum_{l = 0}^\infty b_l \, z^l$ formally to get a new
power series $\sum_{n = 0}^\infty c_n \, z^n$, where
\begin{equation}
\label{cauchy product}
	c_n = \sum_{j = 0}^n a_j \, b_{n - j}.
\end{equation}
This is the \emph{Cauchy product},\index{Cauchy products} and one can
just as well say that the Cauchy product of
\begin{equation}
\label{sum_{j = 0}^infty a_j and sum_{l = 0}^infty b_l}
	\sum_{j = 0}^\infty a_j \quad\hbox{and}\quad\sum_{l = 0}^\infty b_l
\end{equation}
is the series
\begin{equation}
\label{sum_{n = 0}^infty c_n}
	\sum_{n = 0}^\infty c_n.
\end{equation}
It is easy to check that
\begin{equation}
	\sum_{n = 0}^r |c_n| \le \Big(\sum_{j = 0}^r |a_j| \Big)
				\, \Big(\sum_{l = 0}^r |b_l| \Big)
\end{equation}
for all nonnegative integers $r$.  Hence the absolute convergence of
(\ref{sum_{j = 0}^infty a_j and sum_{l = 0}^infty b_l}) implies the
absolute convergence of (\ref{sum_{n = 0}^infty c_n}).  One can show
that the product of the sums in (\ref{sum_{j = 0}^infty a_j and sum_{l
= 0}^infty b_l}) is equal to (\ref{sum_{n = 0}^infty c_n}) in this
case.

	Let $U$ be an open set in the complex plane ${\bf C}$, and let
$f(z)$ be a complex-valued function on $U$ which is differentiable at
some $z = x + y \, i \in U$, $x, y \in {\bf R}$.  This means that
there are complex numbers
\begin{equation}
	\frac{\partial f}{\partial x}(z),
		\quad \frac{\partial f}{\partial y}(z)
\end{equation}
such that
\begin{equation}
\label{f(zeta) = f(z) + ...}
	f(\zeta) = f(z) + \frac{\partial f}{\partial x}(z) \, (\xi - x)
			+ \frac{\partial f}{\partial y}(z) \, (\eta - y)
			+ e_z(\zeta) \, |\zeta - z|,
\end{equation}
where $\zeta = \xi + \eta \, i \in U$, $\xi, \eta \in {\bf R}$, and
\begin{equation}
\label{lim_{zeta to z} e_z(zeta) = 0}
	\lim_{\zeta \to z} e_z(\zeta) = 0.
\end{equation}
If we put
\begin{equation}
	\frac{\partial f}{\partial z}(z)
		= \frac{1}{2} \Big(\frac{\partial f}{\partial x}(z)
				- \frac{\partial f}{\partial y}(z) \, i \Big),
 \quad \frac{\partial f}{\partial \overline{z}}(z)
		= \frac{1}{2} \Big(\frac{\partial f}{\partial x}(z)
				+ \frac{\partial f}{\partial y}(z) \, i \Big),
\end{equation}
then we can rewrite (\ref{f(zeta) = f(z) + ...}) as
\begin{equation}
	f(\zeta) = f(z) + \frac{\partial f}{\partial z}(z) \, (\zeta - z)
			+ \frac{\partial f}{\partial \overline{z}}(z)
				\, (\overline{\zeta} - \overline{z})
			+ e_z(\zeta) \, |\zeta - z|,
\end{equation}
where $e_z(\zeta)$ satisfies (\ref{lim_{zeta to z} e_z(zeta) = 0}) as
before.  We say that $f(z)$ is \emph{holomorphic}\index{holomorphic
functions} on $U$ if it is differentiable everywhere in $U$ and
\begin{equation}
	\frac{\partial f}{\partial \overline{z}}(z) = 0
\end{equation}
for every $z \in U$.  In this event we put
\begin{equation}
	f'(z) = \frac{\partial f}{\partial z}(z),
\end{equation}
the complex derivative of $f$ at $z$.

	Constant functions are trivially holomorphic with derivative
equal to $0$.  The identity function $f(z) = z$ is holomorphic on the
whole complex plane, with derivative equal to $1$.  The product rule
implies that $z^n$ is a holomorphic function on ${\bf C}$ for each
positive integer $n$, with derivative equal to $n \, z^{n - 1}$.  If
the power series (\ref{sum_{n = 0}^infty a_n z^n}) has radius of
convergence $R$, then the series
\begin{equation}
\label{sum_{n = 0}^infty n a_n z^n}
	\sum_{n = 0}^\infty n \, a_n \, z^n
\end{equation}
also has radius of convergence equal to $R$, and the function defined
for $|z| < R$ by (\ref{sum_{n = 0}^infty a_n z^n}) is holomorphic and
its derivative is given by the power series (\ref{sum_{n = 0}^infty n
a_n z^n}).  This follows from standard results about differentiating
power series term by term.

\section{The exponential function}
\label{exponential function}
\setcounter{equation}{0}

	Consider the power series
\begin{equation}
	\sum_{n = 0}^\infty \frac{z^n}{n!},
\end{equation}
where $n!$ or ``$n$ factorial'' is the product of the positive
integers from $1$ to $n$, which is interpreted as being equal to $1$
when $n = 0$.  This series converges absolutely for every complex
number $z$, and the sum is denoted $\exp z$.

	By the remarks in the previous section, the \emph{exponential
function}\index{exponential function} $\exp z$ is a continuous
complex-valued function on ${\bf C}$ which is holomorphic and whose
derivative is equal to itself.  The exponential function is
characterized by these properties and the normalization $\exp 0 = 1$.

	One can show that
\begin{equation}
	\exp (z + w) = (\exp z ) \, (\exp w)
\end{equation}
for all complex numbers $z$, $w$.  More precisely, the binomial
theorem gives
\begin{equation}
	\exp (z + w) = \sum_{n = 0}^\infty \frac{(z + w)^n}{n!}
 		= \sum_{n = 0}^\infty \sum_{l = 0}^n 
			\frac{z^l}{l!} \, \frac{w^{n - l}}{(n - l)!},
\end{equation}
which is the same as the Cauchy product of the series defining $\exp
z$, $\exp w$.  Because the series converge absolutely there is no
problem in identifying the Cauchy product with the product of the
exponentials.

	Because $\exp 0 = 1$, we get that
\begin{equation}
	(\exp z) \, (\exp -z) = 1.
\end{equation}
In particular, $\exp z \ne 0$ for every $z \in {\bf C}$.

	If $x \in {\bf R}$, then $\exp x \in {\bf R}$, since the
coefficients in the series expansion for the exponential function are
real numbers.  Observe that $\exp x \ge 1$ when $x \ge 0$, since all
the terms in the series expansion for the exponential function are
nonnegative real numbers.  Consequently,
\begin{equation}
	\exp x > 0
\end{equation}
for every $x \in {\bf R}$.

	We also have that $\overline{\exp z} = \exp \overline{z}$ for
every $z \in {\bf C}$, because the coefficients of the series
expansion are real.  Hence
\begin{equation}
	|\exp z|^2 = (\exp z) \, (\exp \overline{z})
		= \exp (z + \overline{z}) = \exp (2 \re z)
\end{equation}
and therefore
\begin{equation}
	|\exp z| = \exp \re z.
\end{equation}

	In particular, if $t \in {\bf R}$, then
\begin{equation}
	|\exp (i \, t)| = 1.
\end{equation}
One can show that
\begin{equation}
	\exp (i \, t) = \cos t + i \, \sin t,
\end{equation}
by comparing series expansions.  Alternatively, one can view this in
terms of differential equations, using
\begin{equation}
	\frac{d}{dt} \exp (i \, t) = i \, \exp (i \, t).
\end{equation}

	Geometrically, $\exp (i \, t)$ wraps around the unit circle at
unit speed counterclockwise.  Thus $t \in {\bf R}$ is the same as the
length of the oriented arc from $1$ to $\exp (i \, t)$ on the unit
circle, and $\cos t$, $\sin t$ are the projections of this point on
the unit circle onto the real and imaginary axes, as usual.

	Observe in particular that
\begin{equation}
	\exp (2 \, \pi \, i) = 1,
\end{equation}
because the length of the unit circle is equal to $2 \, \pi$.  It
follows that
\begin{equation}
	\exp (z + 2 \, \pi \, i) = \exp z
\end{equation}
for every $x \in {\bf C}$.

\section{Laurent series}
\label{laurent series}
\setcounter{equation}{0}

	A \emph{Laurent series}\index{Laurent series} has the form
\begin{equation}
\label{sum_{n = -infty}^infty a_n z^n}
	\sum_{n = -\infty}^\infty a_n \, z^n
\end{equation}
for some coefficients $a_n \in {\bf C}$, where $n$ runs through all
integers, positive and negative.  This can be considered as a
combination of the two series
\begin{equation}
\label{sum_{n = 0}^infty a_n z^n, sum_{n = 1}^infty a_{-n} z^{-n}}
	\sum_{n = 0}^\infty a_n \, z^n, 
		\quad \sum_{n = 1}^\infty a_{-n} \, z^{-n},
\end{equation}
which are power series in $z$ and $1/z$, respectively.

	The convergence of the Laurent series (\ref{sum_{n =
-infty}^infty a_n z^n}) can be defined in terms of the convergence of
the two series (\ref{sum_{n = 0}^infty a_n z^n, sum_{n = 1}^infty
a_{-n} z^{-n}}) individually.  The radii of convergence of the two
component power series lead to a maximal open annulus
\begin{equation}
\label{z in {bf C} : r < |z| < R }
	\{z \in {\bf C} : r < |z| < R \}
\end{equation}
on which the Laurent series converges absolutely and uniformly
on compact subsets to a holomorphic function, where
\begin{equation}
	0 \le r \le R \le +\infty.
\end{equation}
The series
\begin{equation}
	\sum_{n = -\infty}^\infty n \, a_n \, z^n
\end{equation}
converges absolutely and uniformly on compact subsets of the same open
annulus, and is equal to the complex derivative of the function
defined by (\ref{sum_{n = -infty}^infty a_n z^n}).

	If $r = 0$ and $R > 0$, then the annulus is actually a
punctured disk when $R < \infty$ and a punctured plane when $R =
+\infty$.  If $a_n = 0$ for every $n < 0$, then the Laurent series is
a power series, and one can include $0$ in the domain of convergence
to get a disk or plane.  If $a_n = 0$ for every $n > 0$, then one can
include $z = \infty$ in the domain of convergence of the series, with
the convention that $1/z = 0$.  If $a_n \ne 0$ for at least one and
only finitely many $n < 0$, then one can interpret the series as
taking the value $\infty$ at $z = 0$.  Similarly, if $a_n \ne 0$ for
at least one and only finitely many $n > 0$, then one can interpret
the series as taking the value $\infty$ at $z = \infty$.

	Suppose that $\sum_{j = -\infty}^\infty a_j$ is absolutely
convergent, i.e., that $\sum_{j = -\infty}^\infty |a_j|$ converges,
which is equivalent to the condition that the partial sums $\sum_{j =
-N}^N |a_j|$, $N \ge 0$, be uniformly bounded.  This implies that the
series (\ref{sum_{n = -infty}^infty a_n z^n}) is absolutely convergent
for every $z \in {\bf C}$ with $|z| = 1$, and that the partial sums
converge uniformly to a continuous complex-valued function on the unit
circle.

	Suppose also that $\sum_{l = -\infty}^\infty b_l$ is an
absolutely convergent doubly-infinite series of complex numbers.  The
\emph{Cauchy product}\index{Cauchy products} of $\sum a_j$, $\sum b_l$
is given by
\begin{equation}
	\sum_{n = -\infty}^\infty c_n,
		\quad c_n = \sum_{j = -\infty}^\infty a_j \, b_{n - j}.
\end{equation}
Using the absolute convergence of $\sum a_j$, $\sum b_l$ one can check
that the series defining $c_n$ converges absolutely for every $n$,
that the series $\sum c_n$ converges absolutely, and that
\begin{equation}
	\sum_{n = -\infty}^\infty |c_n|
		\le \Big(\sum_{j = -\infty}^\infty |a_j| \Big)
			\Big(\sum_{l = -\infty}^\infty |b_l| \Big).
\end{equation}

	Moreover, under these conditions we have that
\begin{equation}
	\sum_{n = -\infty} c_n 
		= \Big(\sum_{j = -\infty}^\infty a_j \Big)
				\Big(\sum_{l = -\infty}^\infty b_l \Big).
\end{equation}
One can show this by approximating the infinite series by finite sums.

	For the same reasons, the product of $\sum a_j \, z^j$, $\sum
b_l \, z^l$ is equal to $\sum c_n \, z^n$ for every $z \in {\bf C}$
with $|z| = 1$.  Analogous statements hold for products of holomorphic
functions defined on an annulur region (\ref{z in {bf C} : r < |z| < R
}) by absolutely convergent Laurent series.

\section{Polynomials}
\label{polynomials}
\setcounter{equation}{0}

	A \emph{holomorphic polynomial} on ${\bf C}$ can be expressed
as
\begin{equation}
	a_n \, z^n + a_{n-1} \, z^{n-1} + \cdots + a_1 \, z + a_0,
\end{equation}
where the coefficients $a_0, \ldots, a_n$ are complex numbers.  Such a
polynomial defines a holomorphic function on the complex plane.

	Let us think of the complex variable $z$ as being $x + y \,
i$, where $x$, $y$ are independent real variables, which amounts to
identifying the complex plane with ${\bf R}^2$.  In general
polynomials on ${\bf R}^2$ are given by finite linear combinations of
the real monomials $x^j \, y^l$, $j, l \ge 0$.  Let us continue to use
complex coefficients here, which ensures that holomorphic polynomials
in $z$ are polynomials in the two real variables $x$, $y$.

	Equivalently, a general polynomial of $x$, $y$ is a finite
linear combination of monomials of the form $z^j \, \overline{z}^l$
with complex coefficients.  The holomorphic polynomials are then just
those for which $\overline{z}$ is not necessary.

	The \emph{Laplace operator}\index{Laplace operator} $\Delta$
on the plane is defined as usual by
\begin{equation}
	\Delta = 
	  \frac{\partial^2}{\partial x^2} + \frac{\partial^2}{\partial y}^2,
\end{equation}
which is the same as
\begin{equation}
	\Delta = 
   4 \frac{\partial}{\partial z} \frac{\partial}{\partial \overline{z}}.
\end{equation}
A twice continuously-differentiable complex-valued function $f$ on an
open set $U$ in the plane is said to be \emph{harmonic}\index{harmonic
functions} if
\begin{equation}
	\Delta f = 0
\end{equation}
on $U$.  Holomorphic functions are continuously differentiable of all
orders, by general results, and automatically harmonic.

	\emph{Conjugate holomorphic functions} on open sets in the
plane are defined in the same way as holomorphic functions but with
the differential equation
\begin{equation}
	\frac{\partial}{\partial \overline{z}} f(z) = 0.
\end{equation}
Equivalently, $f(z)$ is conjugate holomorphic if $\overline{f(z)}$ is
holomorphic, which is also equivalent to $f(\overline{z})$ being
holomorphic on the complex conjugate of the domain of $f$.

	Let us note too that $f(z)$ is holomorphic if and only if
$\overline{f(\overline{z})}$ is holomorphic on the complex conjugate
of the domain of $f$.  At any rate, both holomorphic and conjugate
holomorphic functions are automatically harmonic.

	A \emph{harmonic polynomial} on ${\bf C}$ is a polynomial
which defines a harmonic function on ${\bf C}$.  Because holomorphic
and conjugate holomorphic functions are harmonic, every polynomial
which is a linear combination of $z^j$'s and $\overline{z}^l$'s is
harmonic.  Conversely every harmonic polynomial is of this form.

	A general polynomial on ${\bf C}$ is a linear combination of
monomials of the form $z^j \, \overline{z}^l$, which is the same as a
linear combination of monomials $z^j \, |z|^{2 r}$, $\overline{z}^l \,
|z|^{2 r}$, $j, l, r \ge 0$.  It follows that every polynomial on
${\bf C}$ agrees with a harmonic polynomial on the unit circle.

	Let us write ${\bf T}$\index{#T@${\bf T}$} for the unit
circle, i.e., the set of $z \in {\bf C}$ with $|z| = 1$.  If $f(z)$ is
a continuous complex-valued function on ${\bf T}$, then $f(\exp i \,
t)$ is a continuous function on the real line which is periodic with
period $2 \, \pi$.  Every $2 \, \pi$-periodic continuous function on
the real line corresponds to a continuous function on the unit circle
in this way.

	For a continuous complex-valued function $f(z)$ on ${\bf T}$,
the integral
\begin{equation}
	\int_{\bf T} f(z) \, |dz|
\end{equation}
of $f$ over the unit circle with respect to arc length is the same as
\begin{equation}
	\int_0^{2 \pi} f(\exp (i \, t)) \, dt
\end{equation}
as an ordinary Riemann integral on the real line.  We shall typically
normalize these integrals by dividing by $2 \, \pi$.

	If $f(z) = z^n$ for some $n \ge 1$, then $f(\exp (i \, t)) =
\exp (i \, n \, t)$ can be expressed as $(1/i \, n) (d / dt) \, \exp
(i \, n \, t)$, and therefore
\begin{equation}
	\int_{\bf T} z^n \, |dz| = 0.
\end{equation}
Similarly,
\begin{equation}
	\int_{\bf T} \overline{z}^n \, |dz| = 0
\end{equation}
for every $n \ge 1$.  For $n = 0$ we note that the normalized integral
of the constant function $1$ is equal to $1$.

	The standard integral Hermitian inner product on the vector
space $\mathcal{C}({\bf T})$\index{$C(T)$@$\mathcal{C}({\bf T})$} of
continuous complex-valued functions on the unit circle is defined by
\begin{equation}
	\langle f_1, f_2 \rangle_{\bf T} 
   = \frac{1}{2 \, \pi} \int_{\bf T} f_1(z) \, \overline{f_2(z)} \, |dz|,
\end{equation}
$f_1, f_2 \in \mathcal{C}({\bf T})$.  The remarks of the preceding
paragraph show that the family of functions on ${\bf T}$ consisting of
$z^j$, $j \ge 1$, $\overline{z}^l$, $l \ge 1$, and the constant
function $1$ are orthonormal with respect to this inner product.

	Suppose that
\begin{equation}
	h(z) = a_n \, z^n + \cdots + a_1 \, z + a_0
		+ a_{-1} \, \overline{z} + \cdots + a_{-n} \, \overline{z}^n
\end{equation}
is a harmonic polynomial on ${\bf C}$.  The coefficients of $h$ can be
expressed as
\begin{equation}
  a_j = \frac{1}{2 \, \pi} \, \int_{\bf T} h(z) \, \overline{z}^j \, |dz|
  	 = \frac{1}{2 \, \pi} \, \int_{\bf T} h(z) \, z^{-j} \, |dz|,
\end{equation}
where the two integrals are the same because $\overline{z} = 1/z$ when
$|z| = 1$.  In particular, $h(z)$ is uniquely determined by its
restriction to the unit circle.

\section{Harmonic functions}
\label{harmonic functions}
\setcounter{equation}{0}

	Let $\sum_{n = -\infty}^\infty a_n$ be a doubly-infinite series
of complex numbers, and consider
\begin{equation}
\label{sum_{n = 0}^infty a_n z^n + sum_{n = 1}^infty a_{-n} overline{z}^n}
	\sum_{n = 0}^\infty a_n \, z^n
		+ \sum_{n = 1}^\infty a_{-n} \, \overline{z}^n.
\end{equation}
One can view this as a sum of two power series, one in $z$ and the
other in $\overline{z}$.

	Suppose that these power series have radii of convergence $\ge
1$.  This means that $|a_n| \, r^{|n|}$ is uniformly bounded over all
integers $n$ for every $r > 0$ such that $r < 1$.  Equivalently,
\begin{equation}
	\sum_{n = -\infty}^\infty |a_n| \, t^{|n|}
\end{equation}
converges for every $t > 0$ such that $t < 1$.

	It follows that (\ref{sum_{n = 0}^infty a_n z^n + sum_{n =
1}^infty a_{-n} overline{z}^n}) converges absolutely at every point in
the open unit disk
\begin{equation}
	\{z \in {\bf C} : |z| < 1 \}.
\end{equation}
Moreover, the partial sums of the series converge uniformly on compact
subsets of the open unit disk, and defines a continuous complex-valued
function $h(z)$ there.

	By standard results about power series, $h(z)$ is infinitely
differentiable on the open unit disk.  One can differentiate the
series (\ref{sum_{n = 0}^infty a_n z^n + sum_{n = 1}^infty a_{-n}
overline{z}^n}) term by term, and it follows that $h(z)$ is a harmonic
function on the open unit disk.

	For each $r > 0$ such that $r < 1$ and every integer $\ell$,
\begin{equation}
	a_\ell \, r^{|\ell |}
  = \frac{1}{2 \, \pi} \int_{\bf T} h(r \, z) \, z^{-\ell} \, |dz|
  = \frac{1}{2 \, \pi} \int_{\bf T} h(r \, z) \, \overline{z}^\ell \, |dz|.
\end{equation}
This follows by expanding $h(r \, z)$ into a series and interchanging
the order of integration and summation, which is allowed because of
uniform convergence.  The sum of integrals reduces to a single term as
in the previous section.

	Suppose in addition that $h(z)$ has a continuous extension to
the closed unit disk
\begin{equation}
	\{z \in {\bf C} : |z| \le 1 \}.
\end{equation}
A sufficient condition for this to occur is that the series
\begin{equation}
	\sum_{n = -\infty}^\infty |a_n|
\end{equation}
converges, in which event (\ref{sum_{n = 0}^infty a_n z^n + sum_{n =
1}^infty a_{-n} overline{z}^n}) converges absolutely at every point in
the closed unit disk and the partial sums converge uniformly on the
closed unit disk.  Note that the series (\ref{sum_{n = 0}^infty a_n
z^n + sum_{n = 1}^infty a_{-n} overline{z}^n}) is the same as
(\ref{sum_{n = -infty}^infty a_n z^n}) when $|z| = 1$.  At any rate,
if $h(z)$ has a continuous extension to the closed unit disk, also
denoted $h(z)$, then
\begin{equation}
	a_\ell
	 = \frac{1}{2 \, \pi} \int_{\bf T} h(z) \, z^{-\ell} \, |dz|
	 = \frac{1}{2 \, \pi} \int_{\bf T} h(z) \, \overline{z}^\ell \, |dz|
\end{equation}
for every integer $\ell$.  This follows by sending $r \to 1$ in the
earlier formula, although if the series of coefficients converge
absolutely then one can just apply the same argument directly with $r
= 1$.

	Conversely, let $f(z)$ be a continuous complex-valued function
on the unit circle ${\bf T}$.  For every integer $\ell$, put
\begin{equation}
	a_\ell = \frac{1}{2 \, \pi} \int_{\bf T} f(w) \, w^{-\ell} \, |dw|
  = \frac{1}{2 \, \pi} \int_{\bf T} f(w) \, \overline{w}^\ell \, |dw|.
\end{equation}
Thus
\begin{equation}
	|a_\ell| \le \frac{1}{2 \, \pi} \int_{\bf T} |f(w)| \, |dw|
\end{equation}
for each $\ell$, and hence (\ref{sum_{n = 0}^infty a_n z^n + sum_{n =
1}^infty a_{-n} overline{z}^n}) converges absolutely at every point in
the open unit disk with these choices of coefficients.  Let $h(z)$ be
the value of (\ref{sum_{n = 0}^infty a_n z^n + sum_{n = 1}^infty
a_{-n} overline{z}^n}) at $z$ when $|z| < 1$.  We would like to show
that $h(z)$ extends continuously to the closed unit disk with $h(z) =
f(z)$ when $|z| = 1$.

	For $z, w \in {\bf C}$ with $|z| < 1$ and $|w| = 1$, put
\begin{equation}
	P(z, w) = \frac{1}{2 \, pi} \sum_{n = 0}^\infty z^n \, \overline{w}^n
	  + \frac{1}{2 \, \pi} \sum_{n = 1}^\infty \overline{z}^n \, w^n.
\end{equation}
This is the \emph{Poisson kernel}\index{Poisson kernels} associated to
the unit disk.  By construction,
\begin{equation}
	h(z) = \int_{\bf T} f(w) \, P(z, w) \, |dw|
\end{equation}
when $|z| < 1$.  Observe that
\begin{equation}
	P(z, w) = \frac{1}{2 \, \pi} 
		 + \frac{1}{\pi} \re \sum_{n = 1}^\infty z^n \, \overline{w}^n,
\end{equation}
since $z^n \, \overline{w}^n$ and $\overline{z}^n \, w^n$ are complex
conjugates of each other.  In particular, $P(z, w)$ is actually real
for every $z$, $w$.

	Therefore
\begin{equation}
	P(z, w) = \frac{1}{2 \, \pi} 
    + \frac{1}{\pi} \re \frac{z \, \overline{w}}{1 - z \, \overline{w}},
\end{equation}
by summing the geometric series.  Equivalently,
\begin{equation}
	P(z, w) = 
 \frac{1}{2 \, \pi} \re \frac{1 + z \, \overline{w}}{1 - z \overline{w}}.
\end{equation}
Using
\begin{equation}
	\frac{1}{1 - z \, \overline{w}}
		= \frac{1 - \overline{z} \, w}{|1 - z \, \overline{w}|^2},
\end{equation}
we get that
\begin{equation}
	P(z, w) = \frac{1}{2 \, \pi} \frac{1 - |z|^2}{|z - w|^2}.
\end{equation}
This also employs the fact that $z \, \overline{w} - \overline{z} \,
w$ is imaginary and the restriction to $|w| = 1$.  Consequently,
\begin{equation}
	P(z, w) > 0
\end{equation}
for every $z$, $w$.

	It follows from the definition of the Poisson kernel as a
series that
\begin{equation}
	\int_{\bf T} P(z, w) \, |dw| = 1
\end{equation}
for every $z \in {\bf C}$, $|z| < 1$.  Thus $h(z)$ is an average of
the values of $f$ on the unit circle, weighted according to the
Poisson kernel.  For a fixed $\zeta \in {\bf T}$, $P(z, w)$ is
concentrated near $\zeta$ when $z$ is close to $\zeta$.  For instance,
if $\delta > 0$, $|z - \zeta | < \delta$, and $|w - \zeta | \ge 2 \,
\delta$, then $|z - w| > \delta$ and hence
\begin{equation}
	P(z, w) \le \frac{1}{2 \, \pi \, \delta^2} \, (1 - |z|^2),
\end{equation}
which tends to $0$ as $|z| \to 1$.  One can use this and the
continuity of $f$ to show that
\begin{equation}
	\lim_{z \to \zeta \atop |z| < 1} h(z) = f(\zeta),
\end{equation}
which implies that the extension of $h(z)$ to the closed unit disk
defined by setting $h(\zeta) = f(\zeta)$ when $|\zeta | = 1$ is
continuous, as desired.

\section{Fourier series}
\label{fourier series}
\setcounter{equation}{0}

	Let $f(z)$ be a continuous complex-valued function on the unit
circle ${\bf T}$.  For each integer $\ell$, the $\ell$th \emph{Fourier
coefficient}\index{Fourier coefficients} of $f$ is given by
\begin{equation}
  a_\ell = \frac{1}{2 \, pi} \int_{\bf T} f(w) \, \overline{w}^\ell \, |dw|,
\end{equation}
and the associated \emph{Fourier series}\index{Fourier series} is
\begin{equation}
	\sum_{\ell = -\infty}^\infty a_\ell \, z^\ell,
\end{equation}
$z \in {\bf T}$.

	More precisely, the Fourier series is defined as a formal
series, and convergence of the series is one of the principal
questions.  There are various extensions of Fourier series and
versions of convergence to consider too.

	If $\sum_{n = -\infty}^\infty c_n$ is any series of complex
numbers, then we say that this series is \emph{Abel
summable}\index{Abel summability} if $\sum_{n = -\infty}^\infty |c_n|
\, r^{|n|}$ converges for every $r > 0$ such that $r < 1$ and
\begin{equation}
	\lim_{r \to 1-} \sum_{n = -\infty}^\infty c_n \, r^{|n|}
\end{equation}
exists, in which event the limit is said to be the Abel sum of
$\sum_{n = -\infty}^\infty c_n$.

	It is easy to check that absolutely convergent series are Abel
summable, with the Abel sum equal to the ordinary sum.  One can be a
bit more careful and get a similar conclusion for series which may
converge conditionally.

	For a continuous function $f$ on the unit circle, the Fourier
coefficients satisfy the bound
\begin{equation}
	|a_\ell | \le \frac{1}{2 \, \pi} \int_{\bf T} |f(w)| \, |dw|
\end{equation}
for all $\ell$, and hence $\sum |a_\ell | \, r^{|\ell |}$ converges
absolutely for every $r > 0$ with $r < 1$.  The remarks at the end of
the previous section show that for every $z \in {\bf T}$, the Fourier
series for $f$ at $z$ has Abel sum equal to $f(z)$.

\section{Fourier transforms}
\label{fourier transforms}
\setcounter{equation}{0}

	Let us write $\mathcal{C}({\bf
R})$\index{$C(R)$@$\mathcal{C}({\bf R})$} for the vector space of
continuous complex-valued functions on the real line.

	We say that $f \in \mathcal{C}({\bf R})$ is
\emph{integrable}\index{integrable functions} if
\begin{equation}
\label{int_{bf R} |f(x)| dx < infty}
	\int_{\bf R} |f(x)| \, dx < \infty,
\end{equation}
which is to say that the integrals $\int_a^b |f(x)| \, dx$ of $|f|$
over bounded intervals $[a, b]$ are uniformly bounded over all $a, b
\in {\bf R}$, $a < b$.  The integral in (\ref{int_{bf R} |f(x)| dx <
infty}) can be defined as the supremum of these integrals over bounded
intervals.

	The integrable functions in $\mathcal{C}({\bf R})$ form a
linear subspace denoted $\mathcal{IC}({\bf
R})$.\index{$IC(R)$@$\mathcal{IC}({\bf R})$} If $f \in
\mathcal{IC}({\bf R})$, then the integral $\int_{\bf R} f(x) \, dx$
can be defined as an improper integral, analogous to the way that the
partial sums of an absolutely summable series converge.

	Let us write $\mathcal{C}_{00}({\bf
R})$\index{$C_00(R)$@$\mathcal{C}_{00}({\bf R})$} for the linear
subspace of $\mathcal{C}({\bf R})$ consisting of functions with
\emph{bounded support},\index{bounded support} which is to say
continuous functions $f(x)$ on the real line for which there are $a, b
\in {\bf R}$ such that $f(x) = 0$ when $x \le a$ and when $x \ge b$.
Continuous functions on ${\bf R}$ with bounded support are
automatically integrable, with the integrals reducing to ones on
bounded intervals.

	A continuous function $f(x)$ on the real line is said to be
\emph{bounded}\index{bounded functions} if there is an $A \ge 0$ such
that $|f(x)| \le A$ for every $x \in {\bf R}$.  The bounded functions
in $\mathcal{C}({\bf R})$ form a linear subspace which is denoted
$\mathcal{BC}({\bf R})$\index{$BC(R)$@$\mathcal{BC}({\bf R})$} and
which contains the continuous functions with bounded support.

	Note that the product of a bounded continuous function on
${\bf R}$ and an integrable continuous function on ${\bf R}$ is
integrable.  The product of an arbitrary continuous function on ${\bf
R}$ and a continuous function on ${\bf R}$ that has bounded support
also has bounded support.

	Let $f(x)$ be an integrable continuous complex-valued function
on the real line.  The \emph{Fourier transform}\index{Fourier
transforms} is the function $\widehat{f}(\xi)$ on the real line given
by
\begin{equation}
	\widehat{f}(\xi) = \int_{\bf R} f(x) \, \exp (- \xi \, x \, i) \, dx,
\end{equation}
where the integral makes sense because $\exp (- \xi \, x \, i)$ is a
bounded continuous function on the real line as a function of $x$ for
every $\xi \in {\bf R}$, and hence its product with $f(x)$ is
integrable.  Observe that the Fourier transform of an integrable
function is bounded, and more precisely
\begin{equation}
	|\widehat{f}(\xi)| \le \int_{\bf R} |f(x)| \, dx
\end{equation}
for every $\xi \in {\bf R}$.

	For each $y \in {\bf R}$, put $f_y(x) = f(x - y)$.  Thus $f_y$
is an integrable continuous function on the real line too, and it is
easy to see that its Fourier transform is given by
\begin{equation}
	\widehat{f_y}(\xi) = \widehat{f}(\xi) \, \exp (- \xi \, y \, i).
\end{equation}

\section{Uniform continuity}
\label{uniform continuity}
\setcounter{equation}{0}

	A complex-valued function $f(x)$ on the real line is said to
be \emph{uniformly continuous}\index{uniform continuity} if for each
$\epsilon > 0$ there is a $\delta > 0$ such that $|f(x) - f(y)| <
\epsilon$ for every $x, y \in {\bf R}$ with $|x - y| < \delta$.  There
is a general theorem which states that continuous functions are
uniformly continuous on compact sets, and hence a continuous function
on the real line with bounded support is uniformly continuous.

	Uniformly continuous functions are automatically continuous,
and we let $\mathcal{UC}({\bf R})$\index{$UC(R)$@$\mathcal{UC}({\bf
R})$} be the linear subspace of $\mathcal{C}({\bf R})$ consisting of
uniformly continuous functions.  Let $\mathcal{BUC}({\bf
R})$\index{$BUC(R)$@$\mathcal{BUC}({\bf R})$} be the space of bounded
uniformly continuous complex-valued functions on ${\bf R}$, a linear
subspace of $\mathcal{BC}({\bf R})$.

	The product of two continuous functions on the real line is
also a continuous function, and thus $\mathcal{C}({\bf R})$ is
actually a commutative algebra with respect to pointwise
multiplication.  Clearly $\mathcal{BC}({\bf R})$ is a subalgebra of
$\mathcal{C}({\bf R})$, since the product of two bounded functions is
a bounded function.  One can check that the product of two bounded
uniformly continuous functions is uniformly continuous, which implies
that $\mathcal{BUC}({\bf R})$ is a subalgebra of $\mathcal{BC}({\bf
R})$.

	In general the product of two uniformly continuous functions
is not uniformly continuous.  The sum of two uniformly continuous
functions is uniformly continuous and the product of a uniformly
continuous function with a constant is uniformly continuous, which is
why $\mathcal{UC}({\bf R})$ is a vector space.

	The \emph{supremum norm}\index{supremum norm} of a bounded
continuous function $f$ on the real line is defined by
\begin{equation}
	\|f\|_\infty = \sup \{|f(x)| : x \in {\bf R} \}.
\end{equation}
Observe that
\begin{equation}
	\|f_1 + f_2\|_\infty \le \|f_1\|_\infty + \|f_2\|_\infty
\end{equation}
and
\begin{equation}
	\|f_1 \, f_2\|_\infty \le \|f_1\|_\infty \, \|f_2\|_\infty
\end{equation}
for every $f_1, f_2 \in \mathcal{BC}({\bf R})$.  The \emph{supremum
metric} on $\mathcal{BC}({\bf R})$ is defined by
\begin{equation}
	d_\infty(f_1, f_2) = \|f_1 - f_2\|_\infty
\end{equation}
for $f_1, f_2 \in \mathcal{BC}({\bf R})$.

	It is well known that $\mathcal{BC}({\bf R})$ is
\emph{complete} as a metric space with respect to the supremum metric,
in the sense that every Cauchy sequence in $\mathcal{BC}({\bf R})$
converges.  Also, $\mathcal{BUC}({\bf R})$ is a closed subspace of
$\mathcal{BC}({\bf R})$ with respect to the supremum metric, which
implies that $\mathcal{BUC}({\bf R})$ is complete with respect to the
supremum metric too.

	Let $f$ be a bounded continuous function on the real line, and
for every $y \in {\bf R}$ put $f_y(x) = f(x - y)$, as in the previous
section.  It is easy to see that $f$ is uniformly continuous if and
only if $y \mapsto f_y$ is continuous as a mapping from ${\bf R}$ into
$\mathcal{BC}({\bf R})$ with respect to the supremum metric.  There is
a similar statement for functions which may not be bounded using the
mapping $y \mapsto f_y - f$.

	A complex-valued function $f(x)$ on the real line is said to
be \emph{Lipschitz}\index{Lipschitz functions} if there is a $C \ge 0$
such that
\begin{equation}
	|f(x) - f(y)| \le C \, |x - y|
\end{equation}
for every $x, y \in {\bf R}$.  Lipschitz functions are clearly
uniformly continuous, and the space of Lipschitz functions on ${\bf
R}$ is a linear subspace of $\mathcal{UC}({\bf R})$.  A continuously
differentiable function on the real line is Lipschitz if and only if
its derivative is bounded.

	For every $x \in {\bf R}$, $\exp (- \xi \, x \, i)$ is a
Lipschitz function of $\xi$ on the real line, and more precisely
\begin{equation}
	|\exp (- \xi_1 \, x \, i) - \exp (- \xi_2 \, x \, i)|
		\le |x| \, |\xi_1 - \xi_2|
\end{equation}
for every $\xi_1, \xi_2 \in {\bf R}$.  This follows from the fact that
the derivative of $\exp ( - \xi \, x \, i)$ in $\xi$ is equal to $- i
\, x \, \exp (- \xi \, x \, i)$, which has modulus equal to $|x|$
everywhere.

	Let $f(x)$ be an integrable continuous function on the real
line.  For every $a, b \in {\bf R}$ with $a \le b$ one can check that
\begin{equation}
	\int_a^b f(x) \, \exp (- i \, \xi \, x) \, dx
\end{equation}
is Lipschitz as a function of $\xi$, using the Lipschitz estimates for
the complex exponentials mentioned in the previous paragraph.

	For each $\rho > 0$ there are $a, b \in {\bf R}$ such that $a
\le b$ and
\begin{equation}
	\int_{\bf R} |f(x)| \, dx < \int_a^b |f(x)| \, dx + \rho.
\end{equation}
This implies that
\begin{equation}
	\biggl|\widehat{f}(\xi) 
		- \int_a^b f(x) \, \exp (- \xi \, x \, i) \, dx \biggr| < \rho
\end{equation}
for every $\xi \in {\bf R}$.  It is easy to see from here that
$\widehat{f}(\xi)$ is uniformly continuous.

\section{Inner product spaces}
\label{inner products}
\setcounter{equation}{0}

	Let $V$ be a vector space with complex numbers as scalars.  An
\emph{inner product}\index{inner products} on $V$, or \emph{Hermitian
inner product}, is a complex-valued function $\langle v, w \rangle$
defined for $v, w \in V$ which satisfies the following properties: (i)
$\langle v, w \rangle$ is linear in $v$ for each fixed $w \in V$,
which is to say that
\begin{equation}
	\langle v_1 + v_2, w \rangle
		= \langle v_1, w \rangle + \langle v_2, w \rangle
\end{equation}
for every $v_1, v_2 \in V$, and
\begin{equation}
	\langle \alpha \, v, w \rangle = \alpha \, \langle v, w \rangle
\end{equation}
for every $\alpha \in {\bf C}$ and $v \in V$; (ii) $\langle v, w
\rangle$ is Hermitian symmetric in the sense that
\begin{equation}
	\langle w, v \rangle = \overline{\langle v, w \rangle}
\end{equation}
for every $v, w \in V$; and (iii) $\langle v, w \rangle$ is
positive-definite, which means that
\begin{equation}
	\langle v, v \rangle \ge 0
\end{equation}
for every $v \in V$, and $\langle v, v \rangle = 0$ if and only if $v
= 0$.

	Observe that the Hermitian symmetry of the inner product
implies that $\langle v, v \rangle \in {\bf R}$ for every $v \in V$.

	For example, if $n$ is a positive integer, consider $V = {\bf
C}^n$, the space of $n$-tuples of complex numbers, using
coordinatewise addition and scalar multiplication.  The standard inner
product on ${\bf C}^n$ is defined by
\begin{equation}
	\langle v, w \rangle = \sum_{j=1}^n v_j \, \overline{w_j},
\end{equation}
$v = (v_1, \ldots, v_n)$, $w = (w_1, \ldots, w_n)$.

	Let $V$ be a complex vector space equipped with an inner
product $\langle v, w \rangle$, and for every $v \in V$ put
\begin{equation}
	\|v\| = \langle v, v \rangle^{1/2}.
\end{equation}
If $V = {\bf C}^n$ with the standard inner product, as in the previous
paragraph, then
\begin{equation}
	\|v\| = \Big(\sum_{j=1}^n |v_j|^2 \Big)^{1/2},
\end{equation}
$v = (v_1, \ldots, v_n)$, is the standard norm.

	In general, the \emph{Cauchy--Schwarz
inequality}\index{Cauchy--Schwarz inequality} states that
\begin{equation}
	|\langle v, w \rangle| \le \|v\| \, \|w\|
\end{equation}
for every $v, w \in V$.  One can show this using the fact that
\begin{equation}
	0 \le \langle \alpha \, v + w, \alpha \, v + w \rangle
  = |\alpha |^2 \, \|v\|^2 + 2 \re \alpha \langle v, w \rangle + \|w\|^2
\end{equation}
for every $\alpha \in {\bf C}$.

	By the definition of $\|v\|$,
\begin{equation}
	\|\alpha \, v\| = |\alpha | \, \|v\|
\end{equation}
for every $\alpha \in {\bf C}$ and $v \in V$.  Moreover,
\begin{equation}
	\|v + w \| \le \|v\| + \|w\|
\end{equation}
for every $v, w \in V$, because
\begin{eqnarray}
	\|v + w \|^2 = \langle v + w, v + w \rangle
		& = & \|v\|^2 + 2 \re \langle v, w \rangle + \|w\|^2	\\
		& \le & (\|v\| + \|w\|)^2			\nonumber
\end{eqnarray}
by the Cauchy--Schwarz inequality.

	A pair of vectors $v, w \in V$ are said to be
\emph{orthogonal}\index{orthogonal vectors} if
\begin{equation}
	\langle v, w \rangle = 0,
\end{equation}
in which event we write $v \perp w$.  Observe that $v \perp w$ is
equivalent to
\begin{equation}
	\|v + w\|^2 = \|v\|^2 + \|w\|^2.
\end{equation}

	A collection of vectors $u_1, \ldots u_n \in V$ is said to be
\emph{orthonormal}\index{orthonormal vectors} if $\|u_j\| = 1$ for
each $j$ and $u_j \perp u_l$ for every $j$, $l$ with $j \ne l$.  If
$u_1, \ldots, u_n$ are orthonormal and $\alpha_1, \ldots, \alpha_n$
are complex numbers, then
\begin{equation}
	\biggl\|\sum_{j = 1}^n \alpha_j \, u_j \biggr\|^2
		= \sum_{j = 1}^n |\alpha_j|^2.
\end{equation}

	Let $u_1, \ldots, u_n$ be a collection of orthonormal vectors
in $V$, and let $W$ be the linear subspace of $V$ consisting of linear
combinations of the $u_j$'s.  By definition this means that $W$
consists of the vectors of the form
\begin{equation}
	w = \sum_{j = 1}^n \alpha_j \, u_j
\end{equation}
for some complex numbers $\alpha_1, \ldots, \alpha_n$.  Because of
orthonormality,
\begin{equation}
	\alpha_j = \langle w, u_j \rangle
\end{equation}
for each $j$ in this case.

	For every $v \in V$ put
\begin{equation}
	P(v) = \sum_{j = 1}^n \langle v, u_j \rangle \, u_j.
\end{equation}
Thus $P(v) \in W$ and $P(w) = w$ for every $w \in W$ by the remarks of
the preceding paragraph.

	By construction, $v - P(v) \perp u_j$ for each $j$, and hence
$v - P(v) \perp w$ for every $w \in W$.  In particular, $v - P(v)
\perp P(v)$, and therefore
\begin{equation}
	\|v\|^2 = \|v - P(v)\|^2 + \|P(v)\|^2
		= \|v - P(v)\|^2 + \sum_{j = 1}^n |\langle v, u_j \rangle|^2.
\end{equation}
Consequently,
\begin{equation}
	\sum_{j = 1}^n |\langle v, u_j \rangle|^2 \le \|v\|^2.
\end{equation}

	If $v \in V$ and $w \in W$, then $v - w$ can be expressed as
the sum of $v - P(v)$ and $P(v) - w$, and these two vectors are
orthogonal to each other since $P(v) - w \in W$.  Hence
\begin{equation}
	\|v - w\|^2 = \|v - P(v)\|^2 + \|P(v) - w\|^2,
\end{equation}
which implies that
\begin{equation}
	\|v - P(v)\| \le \|v - w\|
\end{equation}
with equality exactly when $w = P(v)$.

\section{Fourier series, continued}
\label{f-s, continued}
\setcounter{equation}{0}

	As in Section \ref{polynomials},
\begin{equation}
	\langle f_1, f_2 \rangle_{\bf T} 
  = \frac{1}{2 \, \pi} \int_{\bf T} f_1(z) \, \overline{f_2(z)} \, |dz|
\end{equation}
defines an inner product on the vector space $\mathcal{C}({\bf T})$ of
continuous complex-valued functions on ${\bf T}$.  The functions $z^n$
on ${\bf T}$, where $n$ runs through the integers, are orthonormal
with respect to this inner product.

	Let $f$ be a continuous complex-valued function on the unit
circle, and let $a_\ell$ be the $\ell$th Fourier coefficient of $f$.
Equivalently,
\begin{equation}
	a_\ell = \langle f, z^\ell \rangle_{\bf T}.
\end{equation}

	For $N \ge 0$ put
\begin{equation}
	f_N(z) = \sum_{\ell = -N}^N a_\ell \, z^\ell,
\end{equation}
as a continuous function on the unit circle.  Because of
orthonormality,
\begin{equation}
	\frac{1}{2 \, \pi} \int_{\bf T} |f_N(z)|^2 \, |dz|
		= \sum_{\ell = -N}^N |a_\ell|^2.
\end{equation}
Moreover,
\begin{equation}
	\frac{1}{2 \, \pi} \int_{\bf T} |f(z)|^2 \, |dz|
 = \sum_{\ell = -N}^N |a_\ell|^2
	+ \frac{1}{2 \, \pi} \int_{\bf T} |f(z) - f_N(z)|^2 \, |dz|.
\end{equation}

	Thus
\begin{equation}
	\sum_{\ell = -N}^N |a_\ell|^2
		\le \frac{1}{2 \, \pi} \int_{\bf T} |f(z)|^2 \, |dz|
\end{equation}
for each $N$.  Hence
\begin{equation}
	\sum_{\ell = -\infty}^\infty |a_\ell|^2
		\le \frac{1}{2 \, \pi} \int_{\bf T} |f(z)|^2 \, |dz|,
\end{equation}
where the sum of the $|a_\ell|^2$'s converges in particular.

	Consider
\begin{equation}
	\frac{1}{2 \, \pi} \int_{\bf T} |f(z) - f_N(z)|^2 \, |dz|.
\end{equation}
This is the same as the minimum of
\begin{equation}
	\frac{1}{2 \, \pi} \int_{\bf T} |f(z) - p(z)|^2 \, |dz|
\end{equation}
over all functions $p(z)$ on ${\bf T}$ which are linear combinations
of $z^j$, $-N \le j \le N$, with equality exactly when $p(z) =
f_N(z)$.

	Every continuous function on the unit circle is uniformly
continuous, by compactness.  The arguments described at the end of
Section \ref{harmonic functions} imply that the Abel sums\index{Abel
summability} of the Fourier series of $f$ converge to $f$ uniformly on
${\bf T}$.  Using this one can check that $f$ can be uniformly
approximated by finite linear combinations of the $z^j$'s, which is a
version of the Lebesgue -- Weierstrass -- Stone approximation theorem.

	More explicitly, for each $\epsilon > 0$ there is an
$N_\epsilon \ge 0$ and a function $p_\epsilon(z)$ on ${\bf T}$ such
that $p_\epsilon(z)$ is a linear combination of $z^j$, $-N_\epsilon
\le j \le N_\epsilon$, and
\begin{equation}
	|f(z) - p_\epsilon(z)| < \epsilon
\end{equation}
for every $z \in {\bf T}$.  Hence
\begin{equation}
	\frac{1}{2 \, \pi} \int_{\bf T} |f(z) - p_\epsilon(z)|^2 \, |dz| 
		< \epsilon^2.
\end{equation}

	Because of the minimization property,
\begin{equation}
	\frac{1}{2 \, \pi} \int_{\bf T} |f(z) - f_N(z)|^2 \, |dz| < \epsilon^2
\end{equation}
when $N \ge N_\epsilon$.  Therefore
\begin{equation}
 \lim_{N \to \infty} \frac{1}{2 \, \pi} \int_{\bf T} |f(z) - f_N(z)|^2 \, |dz|
	= 0.
\end{equation}
Furthermore,
\begin{equation}
	\frac{1}{2 \, \pi} \int_{\bf T} |f(z)|^2 \, |dz|
		= \sum_{\ell = -\infty}^\infty |a_\ell|^2.
\end{equation}

\section{Almost periodic functions}
\label{almost periodicity}
\setcounter{equation}{0}

	Let $f$ be a bounded continuous complex-valued function on the
real line.  We say that $f$ is \emph{almost periodic}\index{almost
periodic functions} if for each $\epsilon > 0$ there are $y_1, \ldots,
y_n \in {\bf R}$ such that for every $y \in {\bf R}$,
\begin{equation}
	|f(x - y) - f(x - y_j)| < \epsilon
\end{equation}
for some $j$ and all $x \in {\bf R}$.

	Equivalently, $f$ is almost periodic if for each $\epsilon >
0$ there are $h_1, \ldots, h_n$ in $\mathcal{BC}({\bf R})$ such that
for every $y \in {\bf R}$,
\begin{equation}
	|f(x - y) - h_j(x)| < \epsilon
\end{equation}
for some $j$ and all $x \in {\bf R}$.  More precisely, one can check
that this condition for $\epsilon / 2$ implies the previous one for
$\epsilon$.

	In more abstract terms, $f$ is almost periodic if the
collection of translates $f_y(x) = f(x - y)$, $y \in {\bf R}$, is a
precompact set in $\mathcal{BC}({\bf R})$ with respect to the supremum
metric.  Because $\mathcal{BC}({\bf R})$ is complete, a set of bounded
continuous functions is precompact if and only if it is totally
bounded, which means that for each $\epsilon > 0$ the set is contained
in the union of finitely many balls of radius $\epsilon$.  One can
also ask that the balls have their centers contained in the set,
because a ball with radius $\epsilon / 2$ which intersects the set is
contained in a ball with radius $\epsilon$ centered on the set.

	If $f$ is a continuous periodic function on the real line,
then $f$ is uniformly continuous, since any continuous function is
uniformly continuous on compact sets.  Thus $y \mapsto f_y$ is a
continuous periodic mapping from the real line into $\mathcal{BC}({\bf
R})$, whose image is therefore compact, which implies that $f$ is an
almost periodic function.

	The space of almost periodic functions on ${\bf R}$ is denoted
$\mathcal{AP}({\bf R})$.\index{$AP(R)$@$\mathcal{AP}({\bf R})$}
Observe that every translate of an almost periodic function is almost
periodic, since the set of translations would be the same.

	One can check that the sum or product of two almost periodic
functions on the real line is almost periodic.  In particular, the sum
or product of two periodic functions is almost periodic, without
restrictions on the periods.

	Thus $\mathcal{AP}({\bf R})$ is a subalgebra of
$\mathcal{BC}({\bf R})$.  One can check that $\mathcal{AP}({\bf R})$
is closed with respect to the supremum metric, which basically means
that the limit of a sequence of almost periodic functions on the real
line which converges uniformly is almost periodic.

	A continuous function $f(x)$ on ${\bf R}$ is uniformly
continuous if and only if the family of translations $f_y(x)$, $y \in
{\bf R}$, is equicontinuous at $x = 0$, which means that for each
$\epsilon > 0$ there is a $\delta > 0$ such that
\begin{equation}
	|f_y(x) - f_y(0)| < \epsilon
\end{equation}
for every $x, y \in {\bf R}$ with $|x| < \delta$.  Using this
characterization one can check that almost periodic functions are
uniformly continuous.

	Let $\{\xi_j\}_{j=1}^\infty$ be an arbitrary sequence of real
numbers, let $\sum_{j=1}^\infty a_j$ be a sequence of complex numbers
which is absolutely convergent, and consider
\begin{equation}
	\sum_{j = 1}^\infty a_j \exp (\xi_j \, x \, i).
\end{equation}
The partial sums of this series of functions converges uniformly, by
the Weierstrass $M$-test.  It follows that the sum is almost periodic.

	Let $\mathcal{C}_0({\bf
R})$\index{$C_0(R)$@$\mathcal{C}_0({\bf R})$} be the space of
continuous complex-valued functions $f(x)$ on the real line which
``vanish at infinity'' in the sense that for each $\epsilon > 0$ there
is an $R \ge 0$ such that $|f(x)| < \epsilon$ for every $x \in {\bf
R}$ with $|x| \ge R$.  Continuous functions on the real line with
bounded support have this property trivially.

	One can check that the sum and product of two continuous
functions on the real line which vanish at infinity also vanishes at
infinity.  Continuous function on the real line which vanish at
infinity are bounded and uniformly continuous, since continuous
functions are bounded and uniformly continuous on compact sets.

	One can also check that $\mathcal{C}_0({\bf R})$ is closed
with respect to the supremum metric.  In other words, if a sequence of
continuous functions on the real line vanish at infinity and converge
uniformly, then the limit vanishes at infinity too.  Continuous
functions on the real line with bounded support are dense among
continuous functions which vanish at infinity with respect to the
supremum metric, which is to say that $\mathcal{C}_0({\bf R})$ is the
closure of $\mathcal{C}_{00}({\bf R})$ in $\mathcal{BC}({\bf R})$ with
respect to the supremum metric.

	If $f$ is a continuous function on the real line which
vanishes at infinity and is almost periodic, then $f(x) = 0$ for every
$x \in {\bf R}$.  For if $\epsilon > 0$ and $y_1, \ldots, y_n \in {\bf
R}$ have the property that for every $y \in {\bf R}$ there is a $j$
such that $\|f_y - f_{y_j}\|_\infty \le \epsilon$, then one can apply
this with $|y|$ very large to get that for every $r \ge 0$ there is a
$j$ such that $|f(x - y_j)| < 2 \, \epsilon$ when $|x| \le 2 \, r$,
and hence that $|f(x)| < 2 \, \epsilon$ when $|x| \le r$ and $r \ge
\max (|y_1|, \ldots, |y_n|)$.

\section{Step functions}
\label{step functions}
\setcounter{equation}{0}

	By an \emph{interval}\index{intervals} in the real line we
mean a set $I \subseteq {\bf R}$ such that for every $x, z \in I$ and
$y \in {\bf R}$ such that $x < y < z$, we have that $y \in I$.

	This includes open intervals
\begin{equation}
	(a, b) = \{x \in {\bf R} : a < x < b \},
\end{equation}
closed intervals
\begin{equation}
	[a, b] = \{x \in {\bf R} : a \le x \le b \},
\end{equation}
and the half-open, half-closed intervals
\begin{equation}
	(a, b] = \{x \in {\bf R} : a < x \le b \}
\end{equation}
and
\begin{equation}
	[a, b) = \{x \in {\bf R} : a \le x < b \},
\end{equation}
$a, b \in {\bf R}$, $a \le b$.  More precisely, these are bounded
intervals, and there are unbounded intervals which are open and closed
half-lines as well as the real line.

	In general, if $X$ is a set and $E \subseteq X$, then the
\emph{indicator function}\index{indicator functions} associated to $E$
on $X$ is the function ${\bf 1}_E(x)$ which is equal to $1$ when $x
\in E$ and to $0$ when $x \in X \backslash E$.  A \emph{step
function}\index{step functions} on the real line is a function $f(x)$
which can be expressed as
\begin{equation}
\label{f(x) = sum_{j = 1}^n alpha_j {bf 1}_{I_j}(x)}
	f(x) = \sum_{j = 1}^n \alpha_j \, {\bf 1}_{I_j}(x)
\end{equation}
for some complex numbers $\alpha_1, \ldots, \alpha_n$ and intervals
$I_1, \ldots, I_n$.

	Thus step functions are automatically bounded.  If a function
$f(x)$ can be expressed as a sum as in (\ref{f(x) = sum_{j = 1}^n
alpha_j {bf 1}_{I_j}(x)}) where the intervals $I_j$ are bounded, then
we say that $f(x)$ is an \emph{integrable step
function}.\index{integrable functions}

	If $f(x)$ is a step function on the real line and $a, b \in
{\bf R}$, $a \le b$, then the integral
\begin{equation}
	\int_a^b f(x) \, dx
\end{equation}
can be defined in the usual way, since
\begin{equation}
	\int_a^b {\bf 1}_I(x) \, dx = |I \cap [a, b]|
\end{equation}
for any interval $I$, where $|A|$ denotes the length of an interval
$A$.  If $f(x)$ is an integable step function, then
\begin{equation}
	\int_{\bf R} f(x) \, dx
\end{equation}
makes sense, using
\begin{equation}
	\int_{\bf R} {\bf 1}_I(x) \, dx = |I|
\end{equation}
when $I$ is a bounded interval.

	More generally, one can consider piecewise-continuous
functions on ${\bf R}$, in order to have a class of functions which
includes step functions and continuous functions.

	If $a, b \in {\bf R}$ and $a < b$, then
\begin{equation}
	\int_a^b \exp (- \xi \, x \, i) \, dx
	  = \frac{i}{\xi} (\exp (- \xi \, b \, i) - \exp (- \xi \, a \, i))
\end{equation}
when $\xi \ne 0$ and to $b - a$ when $\xi = 0$.  This is a continuous
function on the real line which vanishes at infinity, and in
particular it follows that integrals of families of continuous
periodic functions may not be almost periodic.

	If $f$ is an integrable step function on the real line, then
the Fourier transform of $f$ is defined as usual by
\begin{equation}
	\widehat{f}(\xi) = \int_{\bf R} f(x) \, \exp (- \xi \, x \, i) \, dx.
\end{equation}
This can be expressed as a linear combination of integrals as in the
previous paragraph, and hence $\widehat{f} \in \mathcal{C}_0({\bf
R})$.

	Let $f$ be a continuous integrable function on the real line,
and let $\epsilon > 0$ be given.  Let $a$, $b$ be real numbers such
that $a < b$ and
\begin{equation}
	\int_{\bf R} |f(x)| \, dx < \int_a^b |f(x)| \, dx + \frac{\epsilon}{2}.
\end{equation}
Because $f$ is uniformly continuous on $[a, b]$, there is an
integrable step function $f_1(x)$ on ${\bf R}$ such that
\begin{equation}
	|f_1(x) - f(x)| < \frac{\epsilon}{2 (b - a)}
\end{equation}
when $x \in I$ and $f_1(x) = 0$ when $x \not\in I$, and therefore
\begin{equation}
	\int_{\bf R} |f(x) - f_1(x)| \, dx < \epsilon.
\end{equation}
This implies that
\begin{equation}
	|\widehat{f}(\xi) - \widehat{f_1}(\xi)| < \epsilon
\end{equation}
for every $\xi \in {\bf R}$.  Since $\widehat{f_1} \in
\mathcal{C}_0({\bf R})$ and $\mathcal{C}_0({\bf R})$ is closed with
respect to the supremum metric, it follows that $\widehat{f} \in
\mathcal{C}_0({\bf R})$, which is a version of the Riemann--Lebesgue
lemma.

\section{Norms}
\label{norms}
\setcounter{equation}{0}

	Let $V$ be a complex vector space and let $N(v)$ be a
nonnegative real-valued function on $V$.  We say that $N(v)$ is a
\emph{norm}\index{norms} on $V$ if $N(v) = 0$ if and only if $v = 0$,
$N(v)$ is homogeneous of degree $1$ in the sense that
\begin{equation}
	N(\alpha \, v) = |\alpha| \, N(v)
\end{equation}
for every $\alpha \in {\bf C}$ and $v \in V$, and
\begin{equation}
\label{N(v + w) le N(v) + N(w)}
	N(v + w) \le N(v) + N(w)
\end{equation}
for every $v, w \in V$.

	A set $E \subseteq V$ is said to be \emph{convex}\index{convex
sets} if for every $v, w \in E$ and $t \in {\bf R}$ with $0 \le t \le
1$,
\begin{equation}
	t \, v + (1 - t) \, w \in E.
\end{equation}
If $N(v)$ is a norm on $V$, then
\begin{equation}
\label{v in V : N(v) le 1}
	\{v \in V : N(v) \le 1 \}
\end{equation}
is convex.  Conversely, the triangle inequality (\ref{N(v + w) le N(v)
+ N(w)}) is implied by the combination of the other conditions on $N$
described in the previous paragraph and the convexity of the closed
unit ball (\ref{v in V : N(v) le 1}).

	For example, let $n$ be a positive integer, and suppose that
$V = {\bf C}^n$.  Put
\begin{equation}
	\|v\|_p = \Big(\sum_{j=1}^n |v_j|^p \Big)^{1/p}
\end{equation}
for $p > 0$ and $v \in {\bf C}^n$.  When $p = \infty$ put
\begin{equation}
	\|v\|_\infty = \max (|v_1|, \ldots, |v_n|).
\end{equation}

	Observe that
\begin{equation}
	\|v\|_\infty \le \|v\|_p
\end{equation}
for every $p > 0$ and $v \in {\bf C}^n$.  Since
\begin{equation}
	\|v\|_q^q = \sum_{j=1}^n |v_j|^q
			\le \|v\|_\infty^{q - p} \, \sum_{j=1}^\infty |v_j|^p
			= \|v\|_\infty^{q-p} \, \|v\|_p^p
\end{equation}
when $0 < p \le q < \infty$, we get that
\begin{equation}
	\|v\|_q \le \|v\|_p
\end{equation}
for every $v \in {\bf C}^n$ when $0 < p \le q$.

	It is easy to see directly from the definitions that $\|v\|_p$
defines a norm on ${\bf C}^n$ when $p = 1, \infty$.  When $1 < p <
\infty$ one can use the convexity of the function $r^p$ on the
nonnegative real numbers to show that the closed unit ball associated
to $\|v\|_p$ is convex, and hence that $\|v\|_p$ is a norm.

	If $0 < p < 1$ and $n \ge 2$, then the closed unit ball
associated to $\|v\|_p$ is not convex and $\|v\|_p$ is not a norm.
Using the inequality $(a + b)^p \le a^p + b^p$ for nonnegative real
numbers $a$, $b$ when $0 < p \le 1$, one can check that
\begin{equation}
	\|v + w\|_p^p \le \|v\|_p^p + \|w\|_p^p
\end{equation}
for every $v, w \in {\bf C}^n$, which is a natural alternative to
convexity.  As another alternative, the closed unit ball associated to
$\|v\|_p$ is ``pseudoconvex'' in the sense of several complex
variables, and $\|v\|_p$ is a ``plurisubharmonic'' function instead of
a convex function.

	If $V$ is a complex vector space equipped with an inner
product, then the norm associated to the inner product as in Section
\ref{inner products} is a norm in the general sense considered here.
In particular, $\|v\|_2$ is the norm associated to the standard inner
product on ${\bf C}^n$.

	Now let $V$ be $\mathcal{C}({\bf T})$, the space of continuous
complex-valued functions on the unit circle.  For $1 \le p < \infty$,
put
\begin{equation}
	\|f\|_p = \|f\|_{p, {\bf T}}
	 = \Big(\frac{1}{2 \, \pi} \int_{\bf T} |f(z)|^p \, |dz| \Big)^{1/p}.
\end{equation}
We can extend this to $p = \infty$ using the supremum norm,
\begin{equation}
	\|f\|_\infty = \|f\|_{\infty, {\bf T}}
		= \sup \{|f(z)| : z \in {\bf T} \}.
\end{equation}

	Observe that
\begin{equation}
	\|f\|_p \le \|f\|_\infty
\end{equation}
for every continuous function $f$ on the unit circle.  More generally,
\begin{equation}
	\|f\|_p \le \|f\|_q
\end{equation}
when $1 \le p \le q < \infty$, because of the convexity of the
function $r^{q/p}$ on the nonnegative real numbers.

	It is easy to see directly that $\|f\|_p$ defines a norm on
$\mathcal{C}({\bf T})$ when $p = 1, \infty$.  When $p = 2$, $\|f\|_2$
is the norm associated to the integral inner product considered in
Section \ref{polynomials}.  One can show that $\|f\|_p$ defines a norm
on $\mathcal{C}({\bf T})$ for every $p$, $1 < p < \infty$, by showing
that the associated closed unit ball is a convex set in
$\mathcal{C}({\bf T})$ using the convexity of the function $r^p$ on
the nonnegative real numbers.

	Similarly, if $f$ is a continuous complex-valued function on
the real line with bounded support, put
\begin{equation}
	\|f\|_p = \|f\|_{p, {\bf R}}
		= \Big( \int_{\bf R} |f(x)|^p \, dx \Big)^{1/p}
\end{equation}
when $1 \le p < \infty$, and let $\|f\|_\infty = \|f\|_{\infty, {\bf
R}}$ be the supremum norm of $f$ as in Section \ref{uniform
continuity}.  One can check that $\|f\|_p$ defines a norm on
$\mathcal{C}_{00}(\bf R)$ for every $p$, $1 \le p \le \infty$, and
that
\begin{equation}
	\langle f_1, f_2 \rangle = \langle f_1, f_2 \rangle_{\bf R}
 = \int_{\bf R} f_1(x) \, \overline{f_2(x)} \, dx
\end{equation}
defines an inner product on $\mathcal{C}_{00}({\bf R})$ for which the
associated norm is equal to $\|f\|_2$.

	As in Section \ref{uniform continuity}, the supremum norm can
be defined in the same way on the space $\mathcal{BC}({\bf R})$ of
bounded continuous complex-valued functions on the real line.  The
norm $\|f\|_1$ can be defined in the same way on the space
$\mathcal{IC}({\bf R})$ of integrable continuous complex-valued
functions on ${\bf R}$.

	For $1 \le p < \infty$, let $\mathcal{IC}^p({\bf
R})$\index{$IC^p(R)$@$\mathcal{IC}^p({\bf R})$} be the space of
continuous complex-valued functions $f(x)$ on the real line such that
$|f(x)|^p$ is integrable on ${\bf R}$.  One can check that
$\mathcal{IC}^p({\bf R})$ is a linear subspace of $\mathcal{C}({\bf
R})$, and it clearly contains the continuous functions with bounded
support.  The norm $\|f\|_p$ can be extended to $\mathcal{IC}^p({\bf
R})$ using the same formula as for functions with bounded support.
The product of two square-integrable continuous functions on the real
line is an integrable continuous function on the real line, and the
integral inner product on $\mathcal{C}_{00}({\bf R})$ described above
extends to an inner product on $\mathcal{IC}^2({\bf R})$ by the same
formula, whose associated norm is equal to $\|f\|_2$.

	As another variant of the same notions, let $a$, $b$ be real
numbers with $a < b$, and let $\mathcal{C}([a, b])$\index{$C([a,
b])$@$\mathcal{C}([a, b])$} be the vector space of continuous
complex-valued functions on the interval $[a, b]$.  For every $f \in
\mathcal{C}([a, b])$, put
\begin{equation}
	\|f\|_p = \|f\|_{p, [a, b]}
		= \Big(\int_a^b |f(x)|^p \, dx \Big)^{1/p}
\end{equation}
when $1 \le p < \infty$ and
\begin{equation}
	\|f\|_\infty = \|f\|_{\infty, [a, b]}
		= \sup \{|f(x)| : a \le x \le b \}.
\end{equation}
These are norms on $\mathcal{C}([a, b])$ for the same reasons as in
the previous situations.  Furthermore,
\begin{equation}
	\langle f_1, f_2 \rangle = \langle f_1, f_2 \rangle_{[a, b]}
		= \int_a^b f_1(x) \, \overline{f_2(x)} \, dx
\end{equation}
defines an inner product on $\mathcal{C}([a, b])$ for which the
associated norm is equal to $\|f\|_2$.

\section{Some computations}
\label{some computations}
\setcounter{equation}{0}

	For each $\eta > 0$, put
\begin{equation}
	A_\eta(x) = \exp (- \eta \, |x|).
\end{equation}
Thus $A_\eta$ is a continuous integrable function on the real line,
with
\begin{equation}
	\int_{\bf R} A_\eta(x) \, dx
		= 2 \int_0^\infty \exp ( - \eta \, x) \, dx
			= \frac{2}{\eta}.
\end{equation}

	Observe that
\begin{equation}
	\int_0^\infty \exp (- \eta \, x - \xi \, x \, i) \, dx
		= \frac{1}{\eta + \xi \, i}.
\end{equation}
Similarly,
\begin{equation}
	\int_{-\infty}^0 \exp ( \eta \, x - \xi \, x \, i) \, dx
		= \int_0^\infty \exp (- \eta \, x + \xi \, x \, i) \, dx
		= \frac{1}{\eta - \xi \, i}.
\end{equation}
Hence
\begin{equation}
	\widehat{A_\eta}(\xi)
		= \frac{1}{\eta - \xi \, i} + \frac{1}{\eta + \xi \, i}
		= \frac{2 \, \eta}{\xi^2 + \eta^2}.
\end{equation}

	Using this formula one can check directly that
$\widehat{A_\eta}(\xi)$ is an integrable continuous function on the
real line, and that $\int_{\bf R} \widehat{A_\eta}(\xi) \, d\xi$ does
not depend on $\eta$.  We would like to determine its precise value.

	For every $\zeta \in {\bf C}$ with $\re \zeta > 0$ there is a
unique $\log \zeta \in {\bf C}$ such that
\begin{equation}
	-\pi < \im \log \zeta < \pi
\end{equation}
and
\begin{equation}
	\exp \log \zeta = \zeta.
\end{equation}
Moreover, $\log \zeta$ is a holomorphic function on the right
half-plane and
\begin{equation}
	\frac{\partial}{\partial \zeta} \log \zeta = \frac{1}{\zeta}.
\end{equation}

	If $a, b \in {\bf R}$, $a \le b$,  and $\eta > 0$, then
\begin{eqnarray}
	\int_a^b \frac{1}{\eta + \xi \, i} \, d\xi & = & 
 -i \int_a^b \frac{\partial}{\partial \xi} \log (\eta + \xi \, i) \, d \xi \\
   & = & - i (\log (\eta + b \, i) - \log (\eta + a \, i)).	\nonumber
\end{eqnarray}
Since $\widehat{A_\eta}(\xi) = 2 \re 1 / (\eta + \xi \, i)$,
\begin{equation}
	\int_a^b \widehat{A_\eta}(\xi) \, d\xi
		= 2 \im (\log (\eta + b \, i) - \log (\eta + a \, i)).
\end{equation}
Therefore
\begin{equation}
	\int_{\bf R} \widehat{A_\eta}(\xi) \, d\xi = 2 \, \pi,
\end{equation}
because $\im log (\eta + \xi \, i)$ tends to $- \pi / 2$ as $\xi \to
-\infty$ and to $\pi / 2$ as $\xi \to \infty$.

	This approach to the computation of the integral is quite
pleasant in the way that the integrand is explicitly a derivative.
Alternatively, the computation of the integral is a standard exercise
in basic calculus using trigonometric substitution.

	One can also view the integral geometrically, by projecting
the line $\eta + \xi \, i$, $\xi \in {\bf R}$, in the complex plane
onto the right half of the unit circle, and checking that the
integrand corresponds to twice the element of arclength on the circle.
The integral can be treated as well by applying Cauchy's theorem in
complex analysis to convert the integral on a line to a simpler
integral on a circle.

\section{Fourier inversion}
\label{fourier inversion}
\setcounter{equation}{0}

	Let $f(x)$ be an integrable continuous complex-valued function
on the real line.  We would like to make sense of
\begin{equation}
\label{frac{1}{2 pi} int_{bf R} widehat{f}(xi) exp (xi x i) d xi}
 \frac{1}{2 \, \pi} \int_{\bf R} \widehat{f}(\xi) \, \exp (\xi \, x \, i) 
								\, d \xi
\end{equation}
and show that it is equal to $f(x)$.

	As an analogue of Abel sums,\index{Abel summability} consider
the integrals
\begin{equation}
\label{frac{1}{2 pi} int_{bf R} widehat{f}(xi) A_eta(xi) exp (xi x i) d xi}
	\frac{1}{2 \, \pi} \int_{\bf R} \widehat{f}(\xi) \, A_\eta(\xi)
					\, \exp (\xi \, x \, i) \, d \xi
\end{equation}
for $\eta > 0$.  These integrals make sense, because $A_\eta(\xi)$ is
integrable for each $\eta > 0$ and because $\widehat{f}$ is a bounded
continuous function.  If $\widehat{f}(\xi)$ happens to be integrable,
then one can check that (\ref{frac{1}{2 pi} int_{bf R} widehat{f}(xi)
exp (xi x i) d xi}) is equal to the limit of (\ref{frac{1}{2 pi}
int_{bf R} widehat{f}(xi) A_eta(xi) exp (xi x i) d xi}) as $\eta \to
0$.

	Using the definition of the Fourier transform we can rewrite
(\ref{frac{1}{2 pi} int_{bf R} widehat{f}(xi) A_eta(xi) exp (xi x i) d
xi}) as a double integral
\begin{equation}
	\frac{1}{2 \, \pi} \int_{\bf R} \int_{\bf R} f(y) \, A_\eta(\xi)
				\, \exp (\xi \, (x - y) \, i) \, dy \, d\xi.
\end{equation}
Integrating in $\xi$ first and using the computations from the
previous section, we get that this is equal to
\begin{equation}
\label{int_{bf R} f(y) P_eta (x - y) dy}
	\int_{\bf R} f(y) \, P_\eta (x - y) \, dy,
\end{equation}
where
\begin{equation}
	P_t(y) = \frac{1}{\pi} \frac{t}{y^2 + t^2}
\end{equation}
is the \emph{Poisson kernel}\index{Poisson kernels} associated to the
real line.

	We would like to show that the Poisson integral (\ref{int_{bf
R} f(y) P_eta (x - y) dy}) of $f$ converges to $f(x)$ as $\eta \to 0$.
If we can do this, then it follows that (\ref{frac{1}{2 pi} int_{bf R}
widehat{f}(xi) A_eta(xi) exp (xi x i) d xi}) converges to $f(x)$ as
$\eta \to 0$.

	A key point is that
\begin{equation}
	\int_{\bf R} P_t(y) \, dy = 1
\end{equation}
for each $t > 0$, by the computations in the previous section.  We
would therefore like to show that
\begin{equation}
\label{lim_{eta to 0} int_{bf R} (f(y) - f(x)) P_eta (x - y) dy = 0}
 \lim_{\eta \to 0} \int_{\bf R} (f(y) - f(x)) \, P_\eta (x - y) \, dy = 0.
\end{equation}
The basic idea is that $P_\eta(y - x)$ is concentrated near $x$ and
$f(y) - f(x)$ is small when $y$ is close to $x$ by continuity.

	For each $\delta > 0$,
\begin{equation}
	\lim_{t \to 0} \sup \{P_t(y) : |y| \ge \delta \} = 0.
\end{equation}
One can also check that
\begin{equation}
	\lim_{t \to 0} \int_{|y| \ge \delta } P_t(y) \, dy = 0.
\end{equation}

	Let $\epsilon > 0$ be given, and suppose that $\delta > 0$ has
the property that
\begin{equation}
	|f(y) - f(x)| < \epsilon
\end{equation}
when $|y - x| < \delta$.  It follows that
\begin{equation}
	\int_{|y - x| < \delta} |f(y) - f(x)| \, P_\eta (y - x) \, dy
		< \epsilon \int_{\bf R} P_\eta (y - x) \, dy = \epsilon.
\end{equation}
Using the concentration properties mentioned in the previous
paragraph, one can check that
\begin{equation}
	\lim_{\eta \to 0} \int_{|y - x| \ge \delta} (|f(y)| + |f(x)|) 
			\, P_\eta (y - x) \, dy = 0.
\end{equation}

	Thus (\ref{lim_{eta to 0} int_{bf R} (f(y) - f(x)) P_eta (x -
y) dy = 0}) holds.  Similar arguments show that the Poisson integral
(\ref{int_{bf R} f(y) P_eta (x - y) dy}) of $f$ converges to $f(x)$
uniformly on bounded subsets of the real line, and uniformly
on the whole real line if $f$ is bounded and uniformly continuous.

\section{Several complex variables}
\label{several variables}
\setcounter{equation}{0}

	Fix an integer $n \ge 1$.  By a
\emph{multi-index}\index{multi-indices} we mean an $n$-tuple $\alpha =
(\alpha_1, \ldots, \alpha_n)$ of nonnegative integers.

	For $z = (z_1, \ldots, z_n) \in {\bf C}^n$, put
\begin{equation}
	z^\alpha = z_1^{\alpha_1} \, z_2^{\alpha_2} \cdots z_n^{\alpha_n}.
\end{equation}
We interpret $z_j^{\alpha_j}$ as being equal to $1$ when $\alpha_j =
1$, and hence $z^\alpha = 1$ for every $z \in {\bf C}^n$ when $\alpha
= 0$.

	A \emph{holomorphic polynomial} on ${\bf C}^n$ is given by a
sum
\begin{equation}
	\sum_{\alpha \in A} a_\alpha \, z^\alpha,
\end{equation}
where $A$ is a set of finitely many multi-indices and the $a_\alpha$'s
are complex numbers.

	If $U$ is an open set in ${\bf C}^n$ and $f(z)$ is a
continuously-differentiable complex-valued function on $U$, then the
partial derivatives $(\partial / \partial z_j) f$, $(\partial /
\partial \overline{z}_j) f$ are defined by applying the partial
derivatives in one complex variable described in Section \ref{power
series} to $f$ as a function of $z_j$, $1 \le j \le n$.  We say that
$f(z)$ is \emph{holomorphic}\index{holomorphic functions} on $U$ if
$(\partial / \partial \overline{z}_j) f(z) = 0$ for every $z \in U$
and $1 \le j \le n$.

	Note that $(\partial / \partial z_j ) z_l^p = (p - 1) \,
z_l^{p-1}$ when $j = l$ and $p$ is a positive integer, and is $0$ when
$j \ne l$.  For every $j$, $l$, and $p$, $(\partial / \partial
\overline{z}_j) z_l^p = 0$, and thus holomorphic polynomials are
holomorphic functions on ${\bf C}^n$.

	An infinite sum $\sum_\alpha a_\alpha$ of nonnegative real
numbers using all multi-indices $\alpha$ is said to converge if the
sums $\sum_{\alpha \in E} a_\alpha$ over arbitrary sets $E$ of
finitely many multi-indices $\alpha$ are bounded.  In this event
$\sum_\alpha a_\alpha$ is defined to be the supremum over all such
sums $\sum_{\alpha \in E} a_\alpha$.

	Similarly, an infinite sum $\sum_\alpha a_\alpha$ of complex
numbers is said to converge absolutely if $\sum_\alpha |a_\alpha|$
converges as a sum of nonnegative real numbers.  One can then make
sense of the sum $\sum_\alpha a_\alpha$ by reducing it to a linear
combination of convergent sums of nonnegative real numbers.

	Alternatively, if $\sum_\alpha a_\alpha$ converges absolutely,
then for every $\epsilon > 0$ there is a set $E_\epsilon$ of finitely
many multi-indices $\alpha$ such that
\begin{equation}
	\sum_\alpha |a_\alpha| < \sum_{\alpha \in E_\epsilon} |a_\alpha|
					+ \epsilon.
\end{equation}
If $A$, $B$ are sets of finitely many multi-indices $\alpha$ such that
$E_\epsilon \subseteq A, B$, then
\begin{equation}
 \biggl|\sum_{\alpha \in A} a_\alpha - \sum_{\alpha \in B} a_\alpha \biggr|
	< \epsilon.
\end{equation}
The sum $\sum_\alpha a_\alpha$ can be characterized by the property
that
\begin{equation}
	\biggl|\sum_{\alpha \in A} a_\alpha - \sum_\alpha a_\alpha \biggr|
		\le \epsilon
\end{equation}
for every set $A$ of finitely many multi-indices such that $E_\epsilon
\subseteq A$.

	A \emph{power series}\index{power series} in the $n$ complex
variables $z_1, \ldots, z_n$ is a sum of the form $\sum_\alpha
a_\alpha \, z^\alpha$, where the $a_\alpha$'s are complex numbers and
the sum extends over all multi-indices $\alpha$.  If $\sum_\alpha
|a_\alpha| \, t^\alpha$ converges for some $t = (t_1, \ldots, t_n)$,
where each $t_j$ is a nonnegative real number, then $\sum_\alpha
a_\alpha \, z^\alpha$ converges absolutely for every $z \in {\bf C}^n$
such that $|z_j| \le t_j$ for each $j$, and the sum defines a
continuous function on this set for reasons of uniform convergence.

	Suppose that $r = (r_1, \ldots, r_n)$ is an $n$-tuple of
positive real numbers, and that $\sum_\alpha a_\alpha \, z^\alpha$
converges absolutely for every $z \in {\bf C}^n$ in the open polydisk
where $|z_j| < r_j$ for each $j$.  The sum then defines a smooth
function on the polydisk with the derivatives given by differentiating
the power series term-by-term, again by considerations of uniform
convergence of these power series on compact subsets of the polydisk.
In particular the power series defines a holomorphic function on the
open polydisk.

	If $\sum_\alpha a_\alpha$, $\sum_\beta b_\beta$ are two sums
of complex numbers extending over all multi-indices $\alpha$, $\beta$,
then the \emph{Cauchy product}\index{Cauchy products} of these two
sums is the sum $\sum_\gamma c_\gamma$, where
\begin{equation}
	c_\gamma = \sum_{\alpha + \beta = \gamma} a_\alpha \, b_\beta.
\end{equation}
More precisely, the sum extends over all multi-indices $\alpha$,
$\beta$ such that
\begin{equation}
	\alpha_j + \beta_j = \gamma_j
\end{equation}
for $j = 1, \ldots, n$, and there are only finitely many such
$\alpha$, $\beta$ since the coordinates of these multi-indices are
nonnegative integers.

	If all but finitely many of the $a_\alpha$'s and
$\beta_\gamma$'s are equal to $0$, then the same holds for the
$c_\gamma$'s, and
\begin{equation}
	\sum_\gamma c_\gamma = \Big(\sum_\alpha a_\alpha \Big)
				\, \Big(\sum_\beta b_\beta \Big).
\end{equation}
If the $a_\alpha$'s and $b_\beta$'s are nonnegative real numbers and
$\sum_\alpha a_\alpha$, $\sum_\beta b_\beta$ converge, then the
$c_\gamma$'s are nonnegative real numbers and one can show that
$\sum_\gamma c_\gamma$ converges and is equal to the product of the
sums of the $a_\alpha$'s and $b_\beta$'s.  If $\sum_\alpha a_\alpha$,
$\sum_\beta b_\beta$ are absolutely convergent sums of complex
numbers, then one can show that $\sum_\gamma c_\gamma$ converges
absolutely and is equal to the product of the sums of the $a_\alpha$'s
and $b_\beta$'s.

	The Cauchy product of two power series $\sum_\alpha a_\alpha
\, z^\alpha$, $\sum_\beta b_\beta \, z^\beta$ is equal to $\sum_\gamma
c_\gamma \, z^\gamma$, where the $c_\gamma$'s are obtained from the
Cauchy product of the $a_\alpha$'s and $b_\beta$'s as in the previous
paragraphs.  If $\sum_\alpha a_\alpha \, z^\alpha$, $\sum_\beta
b_\beta \, z^\beta$ converge absolutelym then it follows that
$\sum_\gamma c_\gamma \, z^\gamma$ converges absolutely and is equal
to the product of $\sum_\alpha a_\alpha \, z^\alpha$, $\sum_\beta
b_\beta \, z^\beta$.

	As a basic example, consider the power series $\sum_\alpha
z^\alpha$, where all of the coefficients are equal to $1$.  Formally
this is the same as the product of the $n$ geometric series $\sum_{j_1
= 0}^\infty z_1^{j_1}, \ldots, \sum_{j_n = 0}^\infty z_n^{j_n}$ in the
variables $z_1, \ldots, z_n$.

	In particular, for each $l \ge 0$, the sum of $z^\alpha$ over
all multi-indices $\alpha$ with $0 \le \alpha_j \le l$ for each $j$ is
equal to the product
\begin{equation}
	\Big(\sum_{j_1 = 0}^l z_1^{j_1} \Big) \cdots
		\Big(\sum_{j_n = 0}^l z_n^{j_n} \Big)
\end{equation}
of the corresponding partial sums of the $n$ geometric series in the
variables $z_1, \ldots, z_n$.  Using this one can check that
$\sum_\alpha z^\alpha$ converges absolutely when $|z_j| < 1$ for each
$j$, as a consequence of the absolute convergence of the classical
geometric series $\sum_{j = 0}^\infty a^j$ when $|a| < 1$.  Moreover,
\begin{equation}
	\sum_\alpha z^\alpha = \prod_{j=1}^n \frac{1}{1 - z_j}
\end{equation}
for every $z \in {\bf C}^n$ such that $|z_j| < 1$ for each $j$.

\section{Multiple Fourier series}
\label{multiple f-s}
\setcounter{equation}{0}

	Fix a positive integer $n$, and in this section let us use
arbitrary $n$-tuples of integers as
multi-indices,\index{multi-indices} i.e., elements of ${\bf Z}^n$.
For $z \in {\bf C}^n$ and $\alpha \in {\bf Z}^n$, put
\begin{equation}
	\widetilde{z}^\alpha 
		= \widetilde{z}_1^{\alpha_1} \cdots \widetilde{z}_n^{\alpha_n},
\end{equation}
where $\widetilde{z}_j^{\alpha_j}$ is equal to $z_j^{\alpha_j}$ when
$\alpha_j > 0$, to $\overline{z}_j^{-\alpha_j}$ when $\alpha_j < 0$,
and to $1$ when $\alpha_j = 0$.  Let ${\bf T}^n$ be the
$n$-dimensional torus, which is the set of $z \in {\bf C}^n$ such that
$|z_j| = 1$ for $j = 1, \ldots, n$.  When $z \in {\bf T}^n$,
$\widetilde{z}^\alpha$ is equal to $z^\alpha$, the product of
$z_j^{\alpha_j}$, $1 \le j \le n$, for every $\alpha \in {\bf Z}^n$.

	A general polynomial on ${\bf C}^n$ can be expressed as a sum
of finitely many terms, where each term is a complex multiple of a
product of nonnegative powers of the the real and imaginary parts of
the coordinates of $z = (z_1, \ldots, z_n)$.  Equivalently, a general
polynomial can be expressed as a sum of finitely many complex
multiples of products of nonnegative powers of the $z_j$'s and
$\overline{z}_j$'s.  A twice-continuously differentiable
complex-valued function on an open set in ${\bf C}^n$ is said to be
\emph{polyharmonic}\index{polyharmonic functions} if it is harmonic in
each $z_j$ separately, $1 \le j \le n$.  The \emph{polyharmonic
polynomials} are the polynomials which can be expressed as the sum of
finitely many complex multiples of $\widetilde{z}^\alpha$'s, $\alpha
\in {\bf Z}^n$, and every general polynomial on ${\bf C}^n$ agrees
with a polyharmonic polynomial on ${\bf T}^n$, as one can see by
removing factors of $|z_j|^2$ whenever possible.

	Let $\mathcal{C}({\bf T}^n)$\index{$C(T^n)$@$\mathcal{C}({\bf
T}^n)$} be the vector space of continuous complex-valued functions on
${\bf T}^n$.  If $f_1, f_2 \in \mathcal{C}({\bf T}^n)$, then put
\begin{equation}
	\langle f_1, f_2 \rangle_{{\bf T}^n}
 = \frac{1}{(2 \, \pi)^n} 
		\, \int_{{\bf T}^n} f_1(z) \, \overline{f_2(z)} \, |dz|,
\end{equation}
where the integral over ${\bf T}^n$ is equivalent to an iterated
integral over ${\bf T}$ in each of the $n$ variables.  This is the
standard integral inner product for continuous functions on the
$n$-dimensional torus.  One can check that the restrictions of the
monomials $\widetilde{z}^\alpha$, $\alpha \in {\bf Z}^n$, to ${\bf
T}^n$ are orthonormal with respect to this inner product.

	For every $f \in \mathcal{C}({\bf T}^n)$ and $\alpha \in {\bf
Z}^n$, $a_\alpha = \langle f, z^\alpha \rangle_{{\bf T}^n}$ is the
\emph{Fourier coefficient}\index{Fourier coefficients} of $f$
associated to $\alpha$, and the corresponding \emph{Fourier
series}\index{Fourier series} is given by
\begin{equation}
	\sum_{\alpha \in {\bf Z}^n} a_\alpha \, z^\alpha.
\end{equation}
A priori the Fourier series is a formal expression whose convergence
properties are to be investigated.  Observe that
\begin{equation}
	|a_\alpha| \le \frac{1}{(2 \, \pi)^n} \int_{{\bf T}^n} |f(z)| \, |dz|
\end{equation}
for every $\alpha \in {\bf Z}^n$.  Because of orthonormality of the
$z^\alpha$'s, $\alpha \in {\bf Z}^n$,
\begin{equation}
	\sum_{\alpha \in A} |a_\alpha|^2
		\le \frac{1}{(2 \, \pi)^n} \int_{{\bf T}^n} |f(z)|^2 \, |dz|
\end{equation}
for every set $A$ with finitely many multi-indices $\alpha$.

	By a \emph{polyharmonic power series}\index{polyharmonic power
series} we mean a series of the form $\sum_{\alpha \in {\bf Z}^n}
a_\alpha \, \widetilde{z}^\alpha$ with complex coefficients
$a_\alpha$.  One can think of this as a combination of $2^n$ power
series in the $z_j$'s and their complex conjugates, according to the
signs of the coordinates of $\alpha$.  If the $a_\alpha$'s are
bounded, for instance, then the series converges absolutely for every
$z \in {\bf C}^n$ such that $|z_j| < 1$ for $1 \le j \le n$.  In this
event the series defines a smooth polyharmonic function on the open
unit polydisk, for reasons of uniform convergence on compact
sub-polydisks.

	In particular we can apply this to the Fourier coefficients
$a_\alpha$ of a continuous function $f$ on ${\bf T}^n$.

	For $z$ in the open unit polydisk and $w \in {\bf T}^n$, the
\emph{Poisson kernel}\index{Poisson kernels} $P_n(z, w)$ can be
defined as $1 / (2 \, \pi)^n$ times the sum over $\alpha \in {\bf
Z}^n$ of $\widetilde{z}^\alpha$ times $\widetilde{w}^{-\alpha} =
w^{-\alpha}$.  This is the same as the product of the $1$-dimensional
Poisson kernels $P(z_j, w_j)$, $1 \le j \le n$.  The main point is
that the polyharmonic power series $\sum_{\alpha \in {\bf Z}^n}
a_\alpha \, \widetilde{z}^\alpha$ with coefficients equal to the
Fourier coefficients of $f$ is equal to the integral of $f(w)$ times
$P_n(z, w)$ over $w \in {\bf T}^n$ for each $z$ in the open unit
polydisk.  One can use this to show that as $z$ approaches a point
$\zeta \in {\bf T}^n$, the value of the polyharmonic power series at
$z$ approaches $f(\zeta)$, which is another version of Abel
summability.\index{Abel summability}

	As in Section \ref{f-s, continued}, one can show that
$\sum_{\alpha \in {\bf Z}^n} |a_\alpha|^2$ is equal to $1 / (2 \,
\pi)^n$ times the integral of $|f|^2$ over ${\bf T}^n$, where the sum
is defined to be the supremum of all subsums over finitely many
$\alpha$.  For any set $A$ of finitely many $\alpha \in {\bf Z}^n$,
let $f_A$ be the function on ${\bf T}^n$ which is the sum of $a_\alpha
\, z^\alpha$ over $\alpha \in A$.  One also has that $1 / (2 \,
\pi)^n$ times the integral of $|f - f_A|^2$ over ${\bf T}^n$ is equal
to the sum of $|a_\alpha|^2$ over $\alpha \in {\bf Z}^n \backslash A$,
which is as small as one would like for suitably-large sets $A$.

	Fix a positive integer $n$, and let $\phi$ be a continuous
complex-valued function on the $n$-dimensional torus ${\bf T}^n$.  For
every $a \in {\bf R}^n$ and $z \in {\bf T}^n$, put
\begin{equation}
	f_{a, z}(t) = 
 \phi(z_1 \, \exp (a_1 \, t \, i), \ldots, z_n \, \exp (a_n \, t \, i)).
\end{equation}
By construction, $f_{a, z}$ is a bounded uniformly continuous function
on the real line.  If there is a $b \in {\bf R}$ such that every
coordinate of $a$ is an integer multiple of $b$, then $f_{a, z}$ is a
periodic function on ${\bf R}$.

	In general,
\begin{equation}
	f_{a, z}(t - v) = f_{a, z'}(t),
\end{equation}
where $z'_j = z_j \, \exp (-a_j \, v \, i)$ for every $v \in {\bf R}$.
Using this one can check that $f_{a, z}$ is almost periodic for every
$a \in {\bf R}^n$ and $z \in {\bf T}^n$.  For instance, a sum or
product of $n$ continuous periodic functions on the real line can be
expressed in this way.

\section{Invariant means}
\label{invariant means}
\setcounter{equation}{0}

	Let $f(x)$ be a bounded continuous complex-valued function on
the real line which is almost periodic.  For each $\epsilon > 0$ there
is an $L > 0$ such that
\begin{equation}
	\biggl|\frac{1}{|I|} \int_I f(x - y) \, dx 
		- \frac{1}{|I|} \int_I f(x) \, dx \biggr| < \epsilon
\end{equation}
for every bounded interval $I$ with $|I| \ge L$ and every $y \in {\bf
R}$, since every translate of $f$ can be approximated uniformly by one
of finitely many translates.

	Equivalently, if $I, I' \subseteq {\bf R}$ are bounded
intervals with $|I| = |I'| \ge L$, then
\begin{equation}
\label{|(1/|I|) int_I f(x) dx - (1/|I'|) int_I' f(x) dx| < epsilon}
	\biggl|\frac{1}{|I|} \int_I f(x) \, dx 
		- \frac{1}{|I'|} \int_{I'} f(x) \, dx \biggr| < \epsilon.
\end{equation}
If $I, I' \subseteq {\bf R}$ are bounded intervals such that $|I'| = r
\, |I|$ for some positive integer $r$, then $I'$ is the union of $r$
subintervals $I'_1, \ldots, I'_r$ with length equal to $|I|$,
\begin{equation}
	\frac{1}{|I'|} \int_{I'} f(x) \, dx
 = \frac{1}{r} \, \sum_{p=1}^r \frac{1}{|I'_p|} \int_{I'_p} f(x) \, dx,
\end{equation}
and we can apply the preceding estimate to $I'_p$, $1 \le p \le r$, to
get that (\ref{|(1/|I|) int_I f(x) dx - (1/|I'|) int_I' f(x) dx| <
epsilon}) holds when $|I| \ge L$.

	If $I, I' \subseteq {\bf R}$ are any bounded intervals such
that $|I|, |I'| \ge L$, then
\begin{equation}
	\biggl|\frac{1}{|I|} \int_I f(x) \, dx
		- \frac{1}{|I'|} \int_{I'} f(x) \, dx \biggr| < 3 \, \epsilon,
\end{equation}
since the averages of $f$ over $I$, $I'$ are close to the averages of
$f$ over expanded intervals whose lengths are integer multiples of the
lengths of $I$, $I'$, and the averages of $f$ on these expanded
intervals are approximately the same if the expanded intervals are
sufficiently large and approximately the same.

	To summarize, the averages of $f$ over sufficiently large
intervals are all approximately the same.  It follows that there is an
average $\mu(f)$ in the limit, which is characterized by the property
that for every $\eta > 0$ there is an $L_\eta > 0$ such that
\begin{equation}
	\biggl|\mu(f) - \frac{1}{|I|} \int_I f(x) \, dx \biggr| < \eta
\end{equation}
whenever $I \subseteq {\bf R}$ is a bounded interval which satisfies
$|I| \ge L_\eta$.

	If $f$ is a continuous periodic function on the real line with
period $p$, then the averages of $f$ over arbitrary intervals of
length $p$ are the same.  These averages are also the same as the
averages of $f$ over intervals whose lengths are positive integer
multiples of $p$, and hence their common value is equal to $\mu(f)$.

	If $f$ is a continuous real-valued almost periodic function on
${\bf R}$, then
\begin{equation}
	\mu(f) \in {\bf R}.
\end{equation}
If in addition $f(x) \ge 0$ for every $x \in {\bf R}$, then
\begin{equation}
	\mu(f) \ge 0.
\end{equation}

	It is easy to see that $\mu$ defines a linear mapping from the
space $\mathcal{AP}({\bf R})$ of bounded continuous almost periodic
functions on the real line into the complex numbers.  Furthermore,
\begin{equation}
	|\mu(f)| \le \|f\|_\infty
\end{equation}
for every $f \in \mathcal{AP}({\bf R})$, where $\|f\|_\infty$ denotes
the supremum norm of $f$ on ${\bf R}$.

	If $f \in \mathcal{AP}({\bf R})$ and $h$ is any continuous
complex-valued function on ${\bf C}$, then $h(f(x))$ is almost
periodic too.  One can check this directly from the definitions, using
the fact that $h$ is uniformly continuous on bounded subsets of ${\bf
C}$.  In particular, $|f(x)|$ is almost periodic.

	If $\mu(|f|) = 0$, then $f(x) = 0$ for every $x \in {\bf R}$.
For if $f(x_0) \ne 0$ for some $x_0 \in {\bf R}$, then there are
$\eta, t > 0$ such that
\begin{equation}
	|f(x)| \ge 2 \, \eta
\end{equation}
when $|x - x_0| \le t$.  Because of almost periodicity, every point in
${\bf R}$ is at a bounded distance from some $w \in {\bf R}$
such that
\begin{equation}
	|f(z)| \ge \eta
\end{equation}
when $|z - w| \le t$.  This leads to a positive lower bound for the
averages of $|f|$ on sufficiently large intervals, and hence $\mu(|f|)
> 0$.  It follows that $\mu(|f|^p)^{1/p}$ defines a norm on
$\mathcal{AP}({\bf R})$ when $1 \le p < \infty$.

	Let us focus now on the case where $p = 2$.  Using the
invariant mean $\mu$ we get an inner product
\begin{equation}
	\langle f_1, f_2 \rangle_{\mathcal{AP}({\bf R})}
		= \mu(f_1 \, \overline{f_2})
\end{equation}
on $\mathcal{AP}({\bf R})$ for which the associated norm is
$\mu(|f|^2)^{1/2}$.

	For every $\xi \in {\bf R}$, put
\begin{equation}
	e_\xi(x) = \exp( \xi \, x \, i).
\end{equation}
One can check that $\mu(e_\xi) = 0$ when $\xi \ne 0$.  Therefore the
functions $e_\xi$, $\xi \in {\bf R}$, are orthonormal with respect to
the inner product $\langle f_1, f_2 \rangle_{\mathcal{AP}({\bf R})}$.

	If $f \in \mathcal{A}({\bf R})$ and $\xi_1, \ldots, \xi_n \in
{\bf R}$, then
\begin{equation}
	\sum_{j=1}^n |\langle f, e_{\xi_j} \rangle_{\mathcal{AP}({\bf R})}|^2
		\le \mu(|f|^2),
\end{equation}
because of the orthonormality of the $e_{\xi_j}$'s.  For each
$\epsilon > 0$, the set of $\xi \in {\bf R}$ such that
\begin{equation}
	|\langle f, e_\xi \rangle_{\mathcal{AP}({\bf R})}| \ge \epsilon
\end{equation}
has $\le \mu(|f|^2) / \epsilon^2$ elements, and in particular the set
of $\xi \in {\bf R}$ such that
\begin{equation}
	\langle f, e_\xi \rangle_{\mathcal{AP}({\bf R})} \ne 0
\end{equation}
has only finitely or countably many elements.

\section{Banach spaces}
\label{banach spaces}
\setcounter{equation}{0}

	Let $V$ be a complex vector space equipped with a norm
$\|v\|$.  We say that $V$ is a \emph{Banach space}\index{Banach
spaces} if $V$ is complete as a metric space with the associated
metric $\|v - w\|$.  A complete inner product space is a \emph{Hilbert
space}.\index{Hilbert spaces}

	For each positive integer $n$, ${\bf C}^n$ equipped with the
norm $\|v\|_p$, $1 \le p \le \infty$, as in Section \ref{norms} is a
Banach space, and a Hilbert space when $p = 2$.  The space
$\mathcal{C}({\bf T}^n)$ of countinuous complex-valued functions on
the $n$-dimensional torus ${\bf T}^n$, and the space
$\mathcal{BC}({\bf R})$ of bounded continuous functions on the real
line, are Banach spaces with respect to the supremum norm.

	A closed subspace of a Banach or Hilbert space is a Banach or
Hilbert space, as appropriate.  In particular, the space
$\mathcal{BUC}({\bf R})$ of bounded uniformly continuous functions and
the space $\mathcal{AP}({\bf R})$ of almost periodic functions are
closed linear subspaces of the space $\mathcal{BC}({\bf R})$ of
bounded continuous functions on ${\bf R}$, and hence Banach spaces
with respect to the supremum norm.

	There is a nice characterization of completeness of $V$ in
terms of infinite series.  As usual an infinite series
$\sum_{j=1}^\infty v_j$ with terms $v_j \in V$ for each $j$ is said to
converge if the sequence of partial sums $\sum_{j=1}^n v_j$ converges
in $V$.

	An infinite series $\sum_{j=1}^\infty v_j$ with terms in $V$
is said to converge absolutely if $\sum_{j=1}^\infty \|v_j\|$
converges as an infinite series of nonnegative real numbers.
Equivalently, this holds when the sums $\sum_{j=1}^n \|v_j\|$ are
bounded.

	If $\sum_{j=1}^\infty v_j$ converges absolutely, then for
each $\epsilon > 0$ there is an $L \ge 0$ such that
\begin{equation}
	\sum_{j = l}^n \|v_j\| < \epsilon
\end{equation}
when $n \ge l \ge L$.  Hence
\begin{equation}
	\biggl\|\sum_{j=l}^n v_j \biggr\| < \epsilon
\end{equation}
when $n \ge l \ge L$.

	Thus absolute convergence of $\sum_{j=1}^\infty v_j$ implies
that the sequence of partial sums $\sum_{j=1}^n v_j$ forms a Cauchy
sequence in $V$.  If $V$ is complete, then $\sum_{j=1}^\infty v_j$
converges in $V$.

	Conversely, suppose that every absolutely convergent series in
$V$ converges.  Let $\{v_j\}_{j=1}^\infty$ be a Cauchy sequence in
$V$, which we would like to show converges.

	Because $\{v_j\}_{j=1}^\infty$ is a Cauchy sequence, there is
a subsequence $\{v_{j_l}\}_{l=1}^\infty$ of $\{v_j\}_{j=1}^\infty$
such that
\begin{equation}
	\|v_{j_{l+1}} - v_{j_l}\| \le 2^{-l}
\end{equation}
for each $l$.  Consequently,
\begin{equation}
\label{sum_{l = 1}^infty (v_{j_{l+1}} - v_{j_l})}
	\sum_{l = 1}^\infty (v_{j_{l+1}} - v_{j_l})
\end{equation}
converges absolutely.

	If every absolutely convergent series in $V$ converges, then
(\ref{sum_{l = 1}^infty (v_{j_{l+1}} - v_{j_l})}) converges.  Since
\begin{equation}
	\sum_{l=1}^r (v_{j_{l+1}} - v_{j_l}) = v_{j_{r+1}} - v_{j_1}
\end{equation}
for every positive integer $r$, this means that
$\{v_{j_l}\}_{l=1}^\infty$ converges in $V$.

	A Cauchy sequence with a convergent subsequence also
converges, and thus $\{v_j\}_{j=1}^\infty$ converges in $V$, as
desired.

\section{$\ell^p({\bf Z}^n)$}
\label{l^p(Z^n)}
\setcounter{equation}{0}

	Fix a positive integer $n$.  In the present section it is
again convenient to let multi-indices\index{multi-indices} be
$n$-tuples of arbitrary integers, i.e., elements of ${\bf Z}^n$.

	Let $\mathcal{C}({\bf Z}^n)$\index{$C(Z^n)$@$\mathcal{C}({\bf
Z}^n)$} be the space of families $a = \{a_\alpha\}_{\alpha \in {\bf
Z}^n}$ of complex numbers indexed by ${\bf Z}^n$, which amount to
complex-valued functions on ${\bf Z}^n$.  The support of $a \in
\mathcal{C}({\bf Z}^n)$ is the set of $\alpha \in {\bf Z}^n$ such that
$a_\alpha \ne 0$, and we let $\mathcal{C}_{00}({\bf
Z}^n)$\index{$C_00(Z^n)$@$\mathcal{C}_{00}({\bf Z}^n)$} be the space
of $a \in \mathcal{C}({\bf Z}^n)$ whose support has only finitely many
elements.

	For every $a = \{a_\alpha\}_{\alpha \in {\bf Z}^n} \in
\mathcal{C}_{00}({\bf Z}^n)$, put
\begin{equation}
\label{||a||_p = ||a||_{p, Z^n}	= (sum_{alpha in Z^n} |a_alpha|^p)^{1/p}}
	\|a\|_p = \|a\|_{p, {\bf Z}^n} 
		= \Big(\sum_{\alpha \in {\bf Z}^n} |a_\alpha|^p \Big)^{1/p}
\end{equation}
when $0 < p < \infty$ and
\begin{equation}
\label{||a||_infty = ||a||_{infty, Z^n} = sup {|a_alpha| : alpha in Z^n }}
	\|a\|_\infty = \|a\|_{\infty, {\bf Z}^n}
		= \sup \{|a_\alpha| : \alpha \in {\bf Z}^n \}
\end{equation}
when $p = \infty$.  As in Section \ref{norms}, $\|a\|_p$ defines a
norm on the complex vector space $\mathcal{C}_{00}({\bf Z}^n)$ when $1
\le p \le \infty$, and when $0 < p \le 1$ one has
\begin{equation}
\label{||a + b||_p^p le ||a||_p^p + ||b||_p^p}
	\|a + b\|_p^p \le \|a\|_p^p + \|b\|_p^p,
\end{equation}
$a, b \in \mathcal{C}_{00}({\bf Z}^n)$, as an alternative version of
the triangle inequality.  We also have the inner product
\begin{equation}
\label{langle a, b rangle = langle a, b rangle_{Z^n}, etc}
	\langle a, b \rangle = \langle a, b \rangle_{{\bf Z}^n}
		= \sum_{\alpha \in {\bf Z}^n} a_\alpha \, \overline{b_\alpha}
\end{equation}
for $a, b \in \mathcal{C}_{00}({\bf Z}^n)$, which satisfies $\langle
a, a \rangle = \|a\|_2^2$.

	For $0 < p < \infty$, $\ell^p({\bf
Z}^n)$\index{$l^p(Z^n)$@$\ell^p({\bf Z}^n)$} is defined to be the
space of $a \in \mathcal{C}({\bf Z}^n)$ such that
\begin{equation}
	\sum_{\alpha \in {\bf Z}^n} |c_\alpha|^p < +\infty.
\end{equation}
More precisely, this means that the sums
\begin{equation}
	\sum_{\alpha \in A} |a_\alpha|^p
\end{equation}
over finite sets $A \subseteq {\bf Z}^n$ are bounded, and the sum over
all $\alpha \in {\bf Z}^n$ is defined to be the supremum of these finite sums.

	Let $\ell^\infty({\bf Z}^n)$ be the space of $a \in
\mathcal{C}({\bf Z}^n)$ such that the $a_\alpha$'s are bounded.  We
say that $a \in \mathcal{C}({\bf Z}^n)$ vanishes at infinity if for
each $\epsilon > 0$ there is a finite set $A \subseteq {\bf Z}^n$ such
that
\begin{equation}
	|a_\alpha| < \epsilon
\end{equation}
when $\alpha \not\in A$.  The space of $a \in \mathcal{C}({\bf Z}^n)$
which vanish at infinity is denoted $\mathcal{C}_0({\bf
Z}^n)$,\index{$C_0(Z^n)$@$\mathcal{C}_0({\bf Z}^n)$} and is contained
in $\ell^\infty({\bf Z}^n)$.

	Clearly $\mathcal{C}_{00}({\bf Z}^n)$ is contained in
$\ell^p({\bf Z}^n)$ for every $p$, $0 < p \le \infty$, and in
$\mathcal{C}_0({\bf Z}^n)$.

	Suppose that $a \in \ell^p({\bf Z}^n)$ for some $p$, $0 < p <
\infty$.  For each $\epsilon > 0$, there is a set $A_\epsilon$ of
finitely many $\alpha \in {\bf Z}^n$ such that
\begin{equation}
	\sum_{\alpha \in {\bf Z}^n} |a_\alpha|^p
		< \sum_{\alpha \in A_\epsilon} |a_\alpha|^p + \epsilon^p.
\end{equation}
This implies that $|a_\alpha| < \epsilon$ when $\alpha \not\in
A_\epsilon$, and hence that $a \in \mathcal{C}_0({\bf Z}^n)$.

	Note that (\ref{||a||_p = ||a||_{p, Z^n} = (sum_{alpha in Z^n}
|a_alpha|^p)^{1/p}}), (\ref{||a||_infty = ||a||_{infty, Z^n} = sup
{|a_alpha| : alpha in Z^n }}) carry over to all $a \in \ell^p({\bf
Z}^n)$, $0 < p \le \infty$.  One can check that $\ell^p({\bf Z}^n)$,
$0 < p \le \infty$, and $\mathcal{C}_0({\bf Z}^n)$ are linear
subspaces of $\mathcal{C}({\bf Z}^n)$, and that $\|a\|_p$ defines a
norm on $\ell^p({\bf Z}^n)$ when $1 \le p \le \infty$ and (\ref{||a +
b||_p^p le ||a||_p^p + ||b||_p^p}) holds for $a, b \in \ell^p({\bf
Z}^n)$ when $0 < p \le 1$, by reducing to the earlier inequalities for
finite sums.  If $a \in \ell^p({\bf Z}^n)$, $0 < p < \infty$, then
\begin{equation}
	\|a\|_\infty \le \|a\|_p.
\end{equation}
As in Section \ref{norms},
\begin{equation}
	\Big(\sum_{\alpha \in A} |a_\alpha|^q \Big)^{1/q}
		\le \Big(\sum_{\alpha \in A} |a_\alpha|^p \Big)^{1/p}
\end{equation}
when $a \in \mathcal{C}({\bf Z}^n)$, $0 < p < q < \infty$, and $A
\subseteq {\bf Z}^n$ has only finitely many elements, and it follows
that $\ell^p({\bf Z}^n) \subseteq \ell^q({\bf Z}^n)$ and
\begin{equation}
	\|a\|_q \le \|a\|_p
\end{equation}
when $0 < p < q < \infty$ and $a \in \ell^p({\bf Z}^n)$.

	Let us think of $\ell^p({\bf Z}^n)$ as a metric space with the
metric
\begin{equation}
	\|a - b\|_p
\end{equation}
when $1 \le p \le \infty$ and
\begin{equation}
	\|a - b\|_p^p
\end{equation}
when $0 < p \le 1$.  If a sequence in $\ell^p({\bf Z}^n)$ is a Cauchy
sequence, then the corresponding sequence of values at any $\alpha \in
{\bf Z}^n$ is a Cauchy sequence in ${\bf C}$ and hence converges.  One
can show that the limit is an element of $\ell^p({\bf Z}^n)$ and that
the sequence converges to the limit in $\ell^p({\bf Z}^n)$.  In other
words, $\ell^p({\bf Z}^n)$ is complete for each $p$.

	For $0 < p < \infty$ one can also check that
$\mathcal{C}_{00}({\bf Z}^n)$ is dense in $\ell^p({\bf Z}^n)$.
Similarly, $\mathcal{C}_{00}({\bf Z}^n)$ is dense in
$\mathcal{C}_0({\bf Z}^n)$ with respect to the supremum metric.
Furthermore, $\mathcal{C}_0({\bf Z}^n)$ is equal to the closure of
$\mathcal{C}_{00}({\bf Z}^n)$ in $\ell^\infty({\bf Z}^n)$.  In
particular, $\mathcal{C}_0({\bf Z}^n)$ is complete with respect to the
supremum metric.  There is a more general sense in which
$\mathcal{C}_{00}({\bf Z}^n)$ is dense in $\ell^\infty({\bf Z}^n)$,
i.e., for every $a \in \ell^\infty({\bf Z}^n)$ there is a sequence in
$\mathcal{C}_{00}({\bf Z}^n)$ which is bounded in $\ell^\infty({\bf
Z}^n)$ and which converges to $a$ pointwise on ${\bf Z}^n$.

	If $a \in \ell^1({\bf Z}^n)$, then one can make sense of the
sum
\begin{equation}
\label{sum_{alpha in Z^n} a_alpha}
	\sum_{\alpha \in {\bf Z}^n} a_\alpha.
\end{equation}
When the $a_\alpha$'s are nonnegative real numbers, (\ref{sum_{alpha
in Z^n} a_alpha}) can be defined as the supremum over all finite
subsums.  Every $a \in \ell^1({\bf Z}^n)$ can be expressed as a linear
combination of nonnegative real-valued elements of $\ell^1({\bf
Z}^n)$, which permits one to get (\ref{sum_{alpha in Z^n} a_alpha}) as
a linear combination of sums of nonnegative real numbers.  The
resulting value of (\ref{sum_{alpha in Z^n} a_alpha}) can be
characterized uniquely in terms of approximation by finite sums.  One
can also think of (\ref{sum_{alpha in Z^n} a_alpha}) as a
complex-valued function on $\mathcal{C}_{00}({\bf Z}^n)$ which is
uniformly continuous with respect to the $\ell^1$-metric, and which
therefore has a unique extension to a uniformly continuous function on
all of $\ell^1({\bf Z}^n)$.

	If $a, b \in \ell^2({\bf Z}^n)$, then the product $a \,
\overline{b}$ is an element of $\ell^1({\bf Z}^n)$ and
\begin{equation}
	\sum_{\alpha \in {\bf Z}^n} |a_\alpha| \, |b_\alpha|
		\le \|a\|_2 \, \|b\|_2.
\end{equation}
This follows by applying the Cauchy--Schwarz inequality to sums over
finite subsets of ${\bf Z}^n$.  Hence the inner product (\ref{langle
a, b rangle = langle a, b rangle_{Z^n}, etc}) carries over to all $a,
b \in \ell^2({\bf Z}^n)$.  The norm $\|a\|_2$ is the norm associated
to this inner product.

	More generally, if $0 < p, q, r \le \infty$, $a \in
\ell^p({\bf Z}^n)$, $b \in \ell^q({\bf Z}^n)$, and
\begin{equation}
	\frac{1}{r} = \frac{1}{p} + \frac{1}{q},
\end{equation}
then the product $a \, b \in \ell^r({\bf Z}^n)$, and
\begin{equation}
	\|a \, b\|_r \le \|a\|_p \, \|b\|_q.
\end{equation}
This is a version of \emph{H\"older's inequality},\index{H\"older's
inequality} which is trivial when one of the exponents is infinite, in
which event its reciprocal is interpreted as being equal to $0$.  To
prove it we therefore suppose that $p, q < \infty$, and one can also
make the reductions $r = 1$ and $\|a\|_p = \|b\|_q = 1$.
If $x$, $y$ are nonnegative real numbers, then one can check that
\begin{equation}
	x \, y \le \frac{x^p}{p} + \frac{y^q}{q}.
\end{equation}
This implies that
\begin{equation}
	\sum_{\alpha \in {\bf Z}^n} |a_\alpha| \, |b_\alpha|
		\le \sum_{\alpha \in {\bf Z}^n} \frac{|a_\alpha|^p}{p}
			+ \sum_{\alpha \in {\bf Z}^n} \frac{|b_\alpha|^q}{q}
		= 1,
\end{equation}
as desired.

\section{Measures on ${\bf T}^n$}
\label{measures on T^n}
\setcounter{equation}{0}

	Fix a positive integer $n$.  By a \emph{measure}\index{measures}
on ${\bf T}^n$ we mean a linear function $\mu$ on the vector space
$\mathcal{C}({\bf T}^n)$ of continuous complex-valued functions on
${\bf T}^n$ into the complex numbers which is \emph{bounded} in the
sense that there is an $A \ge 0$ such that
\begin{equation}
\label{|mu(phi)| le A ||phi||_infty}
	|\mu(\phi)| \le A \, \|\phi\|_\infty
\end{equation}
for every $\phi \in \mathcal{C}({\bf T}^n)$, where
\begin{equation}
	\|\phi\|_\infty = \|\phi\|_{\infty, {\bf T}^n}
		= \sup \{|\phi(z)| : z \in {\bf T}^n \}
\end{equation}
is the supremum norm for continuous functions on ${\bf T}^n$.  In this
event we put
\begin{equation}
	\|\mu\|_* = \|\mu\|_{*, {\bf T}^n}
		= \sup \{|\mu(\phi)| : \phi \in \mathcal{C}({\bf T}^n),
					\ \|\phi\|_\infty \le 1 \}.
\end{equation}
Equivalently, (\ref{|mu(phi)| le A ||phi||_infty}) holds with $A =
\|\mu\|_*$, and this is the smallest choice of $A$ with this property.

	It is easy to see that the space of measures on ${\bf T}^n$ is
a vector space with respect to addition and scalar multiplication of
linear functionals on $\mathcal{C}({\bf T}^n)$.  That is to say, if
$\mu_1$, $\mu_2$ are measures on ${\bf T}^n$, then we get another
measure $\mu_1 + \mu_2$ defined by
\begin{equation}
	(\mu_1 + \mu_2)(\phi) = \mu_1(\phi) + \mu_2(\phi)
\end{equation}
for every $\phi \in \mathcal{C}({\bf T}^n)$.  Similarly, if $\mu$ is a
measure on ${\bf T}^n$ and $\alpha \in {\bf C}$, then we get another
measure $c \, \mu$ defined by
\begin{equation}
	(c \, \mu)(\phi) = c \, \mu(\phi)
\end{equation}
for every $\phi \in \mathcal{C}({\bf T}^n)$.  The vector space of
measures on ${\bf T}^n$ is denoted $\mathcal{M}({\bf
T}^n)$.\index{$M(T^n)$@$\mathcal{M}({\bf T}^n)$} One can also check
that $\|\mu\|_*$ defines a norm on $\mathcal{M}({\bf T}^n)$.

	If $f$ is a continuous complex-valued function on ${\bf T}^n$,
then
\begin{equation}
	\mu(\phi) =
 \frac{1}{(2 \, \pi)^n} \int_{{\bf T}^n} \phi(w) \, f(w) \, |dw|
\end{equation}
defines a measure on ${\bf T}^n$, and we have that
\begin{equation}
	\|\mu\|_* = \frac{1}{(2 \, \pi)^n} \int_{{\bf T}^n} |f(w)| \, |dw|.
\end{equation}
One can also allow discontinuous functions $f$ here, as long as one
can make sense of the integrals.  If $z \in {\bf T}^n$, then
\begin{equation}
	\delta_z(\phi) = \phi(z)
\end{equation}
defines a measure on ${\bf T}^n$, the Dirac mass at $z$, and
\begin{equation}
	\|\delta_z\|_* = 1.
\end{equation}

	Suppose that $\{\mu_j\}_{j=1}^\infty$ is a sequence of
measures on ${\bf T}^n$ which is a Cauchy sequence, which means that
for each $\epsilon > 0$ there is an $L \ge 1$ such that
\begin{equation}
	\|\mu_j - \mu_l\|_* < \epsilon
\end{equation}
for every $j, l \ge L$.  It follows that
$\{\mu_j(\phi)\}_{j=1}^\infty$ is a Cauchy sequence in ${\bf C}$ for
every $\phi \in \mathcal{C}({\bf T}^n)$, and hence that
$\{\mu_j(\phi)\}_{j=1}^\infty$ converges as a sequence of complex
numbers.  One can check that the limit defines a measure on ${\bf
T}^n$ and that $\{\mu_j\}_{j=1}^\infty$ converges to this measure in
$\mathcal{M}({\bf T}^n)$.  Therefore $\mathcal{M}({\bf T}^n)$ is
complete and hence a Banach space.

	Let $\mu$ be a measure on ${\bf Z}^n$.  For $\alpha \in {\bf
Z}^n$, the corresponding \emph{Fourier coefficient}\index{Fourier
coefficients} of $\mu$ is given by
\begin{equation}
	a_\alpha = \mu(\psi_\alpha),
\end{equation}
where $\psi_\alpha(w) = w^{-\alpha}$.  Observe that
\begin{equation}
	|a_\alpha| \le \|\mu\|_*
\end{equation}
for every $\alpha \in {\bf Z}^n$.

	Because the $a_\alpha$'s are bounded, the \emph{polyharmonic
power series}\index{polyharmonic power series}
\begin{equation}
	\sum_{\alpha \in {\bf Z}^n} a_\alpha \, \widetilde{z}^\alpha
\end{equation}
converges absolutely for every $z \in {\bf C}^n$ such that $|z_j| <
1$, $1 \le j \le n$.  Here $\widetilde{z}^\alpha =
\widetilde{z}_1^{\alpha_1} \cdots \widetilde{z}_n^{\alpha_n}$, with
$\widetilde{z}_j^{\alpha_j} = z_j^{\alpha_j}$ when $\alpha_j \ge 0$
and $= \overline{z_j}^{-\alpha_j}$ when $\alpha_j \le 0$.

	If $\{\phi_j\}_{j=1}^\infty$ is a sequence of continuous
complex-valued functions on ${\bf T}^n$ which converges uniformly to
the continuous function $\phi$, then
\begin{equation}
	\lim_{j \to \infty} \mu(\phi_j) = \mu(\phi).
\end{equation}
This follows from the boundedness of $\mu$, since
\begin{equation}
	|\mu(\phi_j) - \mu(\phi)| = |\mu(\phi_j - \phi)|
			\le \|\mu\|_* \, \|\phi_j - \phi\|_\infty
\end{equation}
for every $j \ge 1$.

	Let $P_n(z, w)$ be the $n$-dimensional Poisson kernel, as in
Section \ref{multiple f-s}.  For each $z$ in the open unit polydisk,
put
\begin{equation}
	P_{n, z}(w) = (2 \, \pi)^n \, P_n(z, w),
\end{equation}
considered as a continuous function on ${\bf T}^n$.  One can check
that
\begin{equation}
	\sum_{\alpha \in {\bf Z}^n} a_\alpha \, \widetilde{z}^\alpha
		= \mu(P_{n, z})
\end{equation}
for all $z$ in the open unit polydisk.  Basically this is the same as
in the case of Fourier coefficients of continuous functions.
Technical matters of applying $\mu$ to a nice sum and getting the same
answer as applying $\mu$ to the individual terms and then summing can
be handled as in the previous paragraph.

	Fix $r_1, \ldots, r_n \in (0, 1)$, and put
\begin{equation}
	r \circ \zeta = (r_1 \, \zeta_1, \ldots, r_n \, \zeta_n)
\end{equation}
for $\zeta \in {\bf T}^n$.  For every continuous complex-valued
function $f$ on ${\bf T}^n$, consider the expression
\begin{equation}
\label{frac{1}{(2 pi)^n} int_{T^n} mu(P_{n, r circ zeta}) f(zeta) |d zeta |}
 \frac{1}{(2 \, \pi)^n} 
	\int_{{\bf T}^n} \mu(P_{n, r \circ \zeta}) \, f(\zeta) \, |d\zeta |,
\end{equation}
which makes sense because $\mu(P_{n, r \circ \zeta})$ is continuous as
a function of $\zeta$ on ${\bf T}^n$.  The same quantity can be
obtained by applying $\mu$ to
\begin{equation}
	\int_{{\bf T}^n} P_n(r \circ \zeta, w) \, f(\zeta) \, |d\zeta |
\end{equation}
as a function of $w$ on ${\bf T}^n$.  We can rewrite this last
integral as
\begin{equation}
	\int_{{\bf T}^n} P_n(r \circ w, \zeta) \, f(\zeta) \, |d\zeta |,
\end{equation}
since
\begin{equation}
	P_n(r \circ \zeta, w) = P_n(r \circ w, \zeta)
\end{equation}
for every $w, \zeta \in {\bf T}^n$.  These Poisson integrals of $f$
converge uniformly to $f$ on ${\bf T}^n$ as $r_j \to 1$, $1 \le j \le
n$, and hence (\ref{frac{1}{(2 pi)^n} int_{T^n} mu(P_{n, r circ zeta})
f(zeta) |d zeta |}) converges to $\mu(f)$ as $r \to (1, \ldots, 1)$,
which gives a version of \emph{Abel summability}\index{Abel
summability} for the \emph{Fourier series} $\sum_{\alpha \in {\bf
Z}^n} a_\alpha \, z^\alpha$ of $\mu$ in this situation.

\section{Measures on ${\bf R}$}
\label{measures on R}
\setcounter{equation}{0}

	A \emph{measure}\index{measures} on the real line is a linear
mapping $\mu$ from the vector space $\mathcal{C}_{00}({\bf R})$ of
continuous complex-valued functions on ${\bf R}$ with bounded support
into the complex numbers which is bounded in the sense that there is
an $A \ge 0$ such that
\begin{equation}
	|\mu(\phi)| \le A \, \|\phi\|_\infty
\end{equation}
for every $\phi \in \mathcal{C}_{00}({\bf R})$, where
$\|\phi\|_\infty$ is the supremum norm of $\phi$ on ${\bf R}$.  If we
put
\begin{equation}
	\|\mu\|_* = \|\mu\|_{*, {\bf R}}
 		= \sup \{|\mu(\phi)| : \phi \in \mathcal{C}_{00}({\bf R}),
					\ \|\phi\|_\infty \le 1\},
\end{equation}
then the previous inequality holds with $A = \|\mu\|_*$, and this is
the smallest value of $A$ which works.

	The space of measures on ${\bf R}$ is denoted
$\mathcal{M}({\bf R})$\index{$M(R)$@$\mathcal{M}({\bf R})$} and is a
vector space with respect to addition and scalar multiplication of
linear functionals.  It is easy to see that $\|\mu\|_*$ defines a norm
on $\mathcal{M}({\bf R})$.  One can also show that $\mathcal{M}({\bf
R})$ is complete and hence a Banach space.

	If $f$ is a continuous integrable complex-valued function on
${\bf R}$, or an integrable step function, then
\begin{equation}
	\mu(\phi) = \int_{\bf R} f(x) \, \phi(x) \, dx
\end{equation}
defines a measure on ${\bf R}$ with norm given by
\begin{equation}
	\|\mu\|_* = \int_{\bf R} |f(x)| \, dx.
\end{equation}
For every $y \in {\bf R}$, the Dirac mass
\begin{equation}
	\delta_y(\phi) = \phi(y)
\end{equation}
defines a measure on ${\bf R}$ with norm equal to $1$.

	More generally, one might consider linear functionals $\mu$ on
$\mathcal{C}_{00}({\bf R}^n)$ which are bounded when restricted to
functions that are supported in a bounded interval in ${\bf R}$.  In
other words, one would ask that for each bounded interval $I \subseteq
{\bf R}$ there is an $A_I \ge 0$ such that
\begin{equation}
	|\mu(\phi)| \le A_I \, \|\phi\|_\infty
\end{equation}
whenever $\phi$ is a continuous complex-valued function on ${\bf R}$
which satisfies $\phi(x) = 0$ for every $x \in {\bf R} \backslash I$.
For example,
\begin{equation}
	\mu(\phi) = \int_{\bf R} \phi(x) \, dx
\end{equation}
has this property with $A_I = |I|$.

	Let us say that a measure $\mu$ on ${\bf R}$ has \emph{bounded
support}\index{bounded support} if there is a closed and bounded
interval $I \subseteq {\bf R}$ such that $\mu(\phi) = 0$ whenever
$\phi \in \mathcal{C}_{00}({\bf R})$ satisfies $\phi(x) = 0$ for every
$x \in I$.  In this event $\mu(\phi_1) = \mu(\phi_2)$ whenever
$\phi_1, \phi_2 \in \mathcal{C}_{00}({\bf R})$ are equal on $I$, and
$\mu(\phi)$ can be defined for any continuous complex-valued function
$\phi$ on $I$ by extending $\phi$ to a continuous function with
bounded support on ${\bf R}$ and applying $\mu$ to the extension.  One
can choose the extension of $\phi$ in such a way that the supremum
norm of the extension is equal to the supremum norm of $\phi$ on $I$.
In particular, $|\mu(\phi)|$ is less than or equal to $\|\mu\|_*$
times the supremum norm of $\phi$ on $I$.

	It turns out that the measures on ${\bf R}$ with bounded
support are dense in $\mathcal{M}({\bf R})$.

	To see this, let $\mu \in \mathcal{M}({\bf R})$ and $\epsilon
> 0$ be given.  Suppose that $\psi \in \mathcal{C}_{00}({\bf R})$,
$\|\psi\|_* \le 1$, and
\begin{equation}
	\|\mu\|_* < |\mu(\psi)| + \epsilon.
\end{equation}
Let $\rho$ be a real-valued continuous function on the real line with
bounded support such that
\begin{equation}
	0 \le \rho(x) \le 1
\end{equation}
for every $x \in {\bf R}$ and $\rho(x) = 1$ when $\psi(x) \ne 0$.  One
can show that
\begin{equation}
	\mu'(\phi) = \mu(\rho \, \phi)
\end{equation}
defines a measure on ${\bf R}$ with bounded support and $\|\mu -
\mu'\|_* < \epsilon$.

	Let us say that a measure $\mu$ on ${\bf R}$ has a
\emph{regular extension}\index{regular extensions} if there is an
extension of $\mu$ as a linear functional to the vector space of
bounded continuous functions on the real line, also denoted $\mu$,
such that
\begin{equation}
	\lim_{j \to \infty} \mu(\phi_j) = \mu(\phi)
\end{equation}
whenever $\{\phi_j\}_{j=1}^\infty$ is a uniformly bounded sequence of
continuous complex-valued functions on the real line which converges
uniformly on bounded intervals to the continuous function $\phi$.
Such an extension would be unique, because for every bounded
continuous complex-valued function $\phi$ on ${\bf R}$ there is a
uniformly bounded sequence $\{\phi_j\}_{j=1}^\infty$ of continuous
functions on ${\bf R}$ such that each $\phi_j$ has bounded support and
the sequence converges uniformly to $\phi$ on every bounded interval
in ${\bf R}$.  By choosing $\phi_j$'s with $\|\phi_j\|_* \le
\|phi\|_*$ for each $j$ we also get that
\begin{equation}
	|\mu(\phi)| \le \|\mu\|_* \, \|\phi\|_*
\end{equation}
for every bounded continuous function $\phi$ on the real line.

	If a measure $\mu$ on ${\bf R}$ has bounded support, then
$\mu$ has a regular extension, as a consequence of the earlier
remarks.  Suppose that $\{\mu_j\}_{j=1}^\infty$ is a sequence of
measures on ${\bf R}$ with regular extensions which converges in the
norm on $\mathcal{M}({\bf R})$ to a measure $\mu$.  In this case one
can check that $\mu$ has a regular extension, obtained by taking the
limit of $\{\mu_j(\phi)\}_{j=1}^\infty$ for every bounded continuous
complex-valued function $\phi$ on ${\bf R}$.  It follows that every
measure on ${\bf R}$ has a regular extension.

	Let $\mu$ be a measure on ${\bf R}$, which has a regular
extension to bounded continuous functions on ${\bf R}$, as in the
previous paragraph, also denoted $\mu$.  For every $\xi \in {\bf R}$,
\begin{equation}
	e_\xi(x) = \exp (\xi \, x \, i)
\end{equation}
is a bounded continuous complex-valued function on ${\bf R}$.  The
\emph{Fourier transform}\index{Fourier transforms} of $\mu$ is the
function $\widehat{\mu}$ on ${\bf R}$ given by
\begin{equation}
	\widehat{\mu}(\xi) = \mu(e_{-\xi}),
\end{equation}
which satisfies
\begin{equation}
	|\widehat{\mu}(\xi)| \le \|\mu\|_*
\end{equation}
for every $\xi \in {\bf R}$.  If $\mu$ has bounded support, then one
can check that $\widehat{\mu}$ is a Lipschitz function on the real
line.  Because the measures with bounded support are dense in
$\mathcal{M}({\bf R})$, it follows that $\widehat{\mu}$ is uniformly
continuous for every measure $\mu$ on ${\bf R}$.

	As in Section \ref{fourier inversion}, consider
\begin{equation}
 \frac{1}{2 \, \pi} \int_{\bf R} \widehat{\mu}(\xi) \, A_\eta(\xi)
				\, \exp (\xi \, x \, i) \, d\xi
\end{equation}
as a version of Abel sums\index{Abel summability} for the inverse
Fourier transform of $\widehat{\mu}$, where $A_\eta(xi) = \exp( - \eta
\, |\xi|)$, $\eta > 0$.  If $P_t(y) = (1 / \pi) \, t / (y^2 + t^2)$ is
the Poisson kernel associated to the real line, then this integral is
equal to $\mu(p_{\eta, x})$, $p_{\eta, x}(y) = P_\eta(x - y)$.  One
can show that this is an integrable continuous function on the real
line which satisfies
\begin{equation}
	\int_{\bf R} \mu(p_{\eta, x}) \, \phi(x) \, dx = \mu(P_\eta(\phi))
\end{equation}
for every bounded continuous function $\phi$ on ${\bf R}$, where
$P_t(\phi)$ is the Poisson integral of $\phi$,
\begin{equation}
	P_t(\phi)(w) = \int_{\bf R} P_t(w - x) \, \phi(w) \, dw.
\end{equation}
This amounts to interchanging the order of integration.  Consequently,
\begin{equation}
 \lim_{\eta \to 0} \int_{\bf R} \mu(p_{\eta, x}) \, \phi(x) \, dx = \mu(\phi),
\end{equation}
because $P_t(\phi)$ is a uniformly bounded family of continuous
functions on the real line when $\phi$ is a bounded continuous
function on ${\bf R}$ which converges to $\phi$ uniformly on bounded
intervals as $t \to 0$.

\section{Convolutions on ${\bf T}^n$}
\label{convolutions, t^n}
\setcounter{equation}{0}

	Fix a positive integer $n$, and let $f(z)$, $g(z)$ be
continuous complex-valued functions on the $n$-dimensional torus.  The
\emph{convolution}\index{convolutions} of $f$, $g$ is the function on
${\bf T}^n$ defined by
\begin{equation}
	(f * g)(z) = \frac{1}{(2 \, \pi)^n}
		\int_{{\bf T}^n} f(w) \, g(z \circ \overline{w}) \, |dw|.
\end{equation}
Here
\begin{equation}
	u \circ v = (u_1 \, v_1, \ldots, u_n \, v_n)
\end{equation}
and
\begin{equation}
	\overline{w} = (\overline{w_1}, \ldots, \overline{w_n})
\end{equation}
for $u, v, w \in {\bf T}^n$.

	Observe that
\begin{equation}
	f * g = g * f,
\end{equation}
as one can check using a change of variables.  Because the torus is
compact, $f$, $g$ are uniformly continuous, and this implies that the
convolution $f * g$ is a continuous function on ${\bf T}^n$.

	Let $a_\alpha$, $b_\alpha$, and $c_\alpha$ be the Fourier
coefficients of $f$, $g$, and $f * g$.  A key feature of the
convolution is the identity
\begin{equation}
	c_\alpha = a_\alpha \, b_\alpha
\end{equation}
for every $\alpha \in {\bf Z}^n$.  Basically, $c_\alpha$ can be
expressed as a double integral which reduces to a product of integrals
after a change of variables.

	For every $z \in {\bf T}^n$,
\begin{equation}
	|f * g(z)| \le \frac{1}{(2 \, \pi)^n}
		\int_{{\bf T}^n} |f(w)| \, |g(z \circ \overline{w})| \, |dw|.
\end{equation}
One can ckeck that
\begin{equation}
	\|f * g\|_1 \le \|f\|_1 \, \|g\|_1,
\end{equation}
by converting a double integral into a product of integrals again, and
that
\begin{equation}
 \|f * g\|_\infty \le \|f\|_1 \, \|g\|_\infty, \ \|f\|_\infty \, \|g\|_1.
\end{equation}

	Let $\mu$ be a measure and $f$ be a continuous complex-valued
function on ${\bf T}^n$.  For every $z \in {\bf T}^n$, put
\begin{equation}
	f_z(w) = f(z \circ \overline{w}).
\end{equation}
The \emph{convolution}\index{convolutions} of $\mu$ and $f$ is the
function on ${\bf T}^n$ defined by
\begin{equation}
	\mu * f(z) = \mu(f_z).
\end{equation}

	Because $f$ is uniformly continuous on ${\bf T}^n$, $\mu * f$
is a continuous function.  We also have the estimate
\begin{equation}
	\|\mu * f\|_\infty \le \|\mu\|_* \, \|f\|_\infty.
\end{equation}

	More generally, suppose that $\mu$, $\nu$ are measures on
${\bf T}^n$.  We would like to begin by defining a measure $\mu \times
\nu$ on ${\bf T}^{2n} \cong {\bf T}^n \times {\bf T}^n$.  Let $f(z,
w)$ be a continuous function on ${\bf T}^n \times {\bf T}^n$, and let
us say how to evaluate $(\mu \times \nu)(f)$.

	For every $z \in {\bf T}^n$, we can apply $\nu$ to $f(z, w)$
as a function of $w$, and get a function of $z$ which is continuous
because of uniform continuity.  If we apply $\mu$ to the resulting
function of $z$, then we get the first definition of $(\mu \times
\nu)(f)$.

	Alternatively, for every $w \in {\bf T}^n$ we can apply $\mu$
to $f(z, w)$ as a function of $z$, and then apply $\nu$ to the
resulting function of $w$.  We would like to show that these two
definitions of $(\mu \times \nu)(f)$ are the same.

	If $f(z, w) = f_1(z) \, f_2(w)$ for continuous functions
$f_1$, $f_2$ on ${\bf T}^n$, then both definitions yield
\begin{equation}
\label{(mu times nu)(f) = mu(f_1) nu(f_2)}
	(\mu \times \nu)(f) = \mu(f_1) \, \nu(f_2).
\end{equation}
By linearity, both definitions agree for sums of products of functions
of $z$, $w$, separately.  Every continuous function on ${\bf T}^n
\times {\bf T}^n$ can be uniformly approximated by such sums, and one
can use this to show that the two definitions of $(\mu \times \nu)(f)$
agree for every continuous function $f(z, w)$.

	Both definitions imply that
\begin{equation}
	|(\mu \times \nu)(f)| 
 		\le \|\mu\|_{*, {\bf T}^n} \, \|\nu\|_{*, {\bf T}^n}
			\, \|f\|_{\infty, {\bf T}^{2n}},
\end{equation}
and hence that $\mu \times \nu$ is a measure on ${\bf T}^{2n}$ with
norm less than or equal to the product of the norms of $\mu$ and $\nu$
on ${\bf T}^n$.  In fact we have that
\begin{equation}
	\|\mu \times \nu\|_{*, {\bf T}^{2n}}
		= \|\mu\|_{*, {\bf T}^n} \, \|\nu\|_{*, {\bf T}^n},
\end{equation}
because of (\ref{(mu times nu)(f) = mu(f_1) nu(f_2)}).

	The \emph{convolution}\index{convolutions} of two measures
$\mu$, $\nu$ on ${\bf T}^n$ is the measure $\mu * \nu$ defined by
applying $\mu \times \nu$ to $f(z \circ w)$ for every continuous
function $f$ on ${\bf T}^n$.  It follows from the previous discussion
that
\begin{equation}
	|(\mu * \nu)(f)| \le \|\mu\|_* \, \|\nu\|_* \, \|f\|_\infty
\end{equation}
for every $f \in \mathcal{C}({\bf T}^n)$, which is to say that $\mu *
\nu$ is bounded and that
\begin{equation}
	\|\mu * \nu\|_* \le \|\mu\|_* \, \|\nu\|_*.
\end{equation}

	If $\mu$ or $\nu$ are given by $(2 \, \pi)^{-n}$ times
integration with a continuous density, then the convolution $\mu *
\nu$ is given by $(2 \, \pi)^{-n}$ times integration with a continuous
density, where the density corresponds to one of the previous
definitions of the convolution, as appropriate.

	For every $\alpha \in {\bf Z}^n$, one can check that the
$\alpha$th Fourier coefficient of $\mu \times \nu$ is equal to the
product of the $\alpha$th Fourier coefficients of $\mu$, $\nu$.

	The convolution of the Dirac mass at $(1, 1, \ldots, 1)$ with
any function or measure is equal to that function or measure.

\section{Convolutions on ${\bf R}$}
\label{convolutions, r}
\setcounter{equation}{0}

	If $f(x)$, $g(x)$ are continuous complex-valued functions on
the real line, then, under suitable additional conditions, the
\emph{convolution}\index{convolutions} of $f$, $g$ is defined by
\begin{equation}
	f * g(x) = \int_{\bf R} f(y) \, g(x - y) \, dx.
\end{equation}

	For instance, this makes sense if $f$, $g$ have bounded
support, in which event $f * g$ has bounded support.  Uniform
continuity of $f$, $g$ imply that $f * g$ is a continuous function on
the real line.

	If at least one of $f$, $g$ has bounded support and the other
is an arbitrary continuous function on ${\bf R}$, then the convolution
$f * g$ is defined as a function on ${\bf R}$.  One can also show that
$f * g$ is continuous in this case.

	If $\mu$ is a measure and $h$ is a bounded continuous function
on ${\bf R}$, then the convolution $\mu * h$ is defined as a function
on ${\bf R}$ by the formula
\begin{equation}
	(\mu * h)(x) = \mu(h_x), \quad h_x(u) = h(x - u).
\end{equation}
This is the same as the previous formula when $\mu$ is given by
integration with a continuous integrable density.  One can check that
$\mu * h$ is continuous using the continuity properties of the regular
extension of $\mu$ to bounded continuous functions on ${\bf R}$.  If
$h$ is bounded and uniformly continuous, then $\mu * h$ is uniformly
continuous.  At any rate, for a bounded continuous function $h$ on
${\bf R}$,
\begin{equation}
	\|\mu * h \|_\infty \le \|\mu\|_* \, \|h\|_\infty.
\end{equation}

	Now suppose that $\mu$, $\nu$ are measures on ${\bf R}$.  If
$f(x, y)$ is a continuous complex-valued function on ${\bf R} \times
{\bf R} \cong {\bf R}^2$, then there are two ways to try to make sense
of $(\mu \times \nu)(f)$, by applying $\mu$ or $nu$ to $f$ as a
function of $x$ or $y$ and then to the other variable.  Because $f$
has bounded support in the plane, the functions on ${\bf R}$ obtained
by applying $\mu$ or $\nu$ to $f(x, y)$ in $x$ or $y$ each have
bounded support.  They are also continuous, for the usual reasons of
uniform continuity.  Hence one can apply $\mu$ or $\nu$ to the
resulting function of one variable.

	After applying $\mu$ or $\nu$ to $f(x, y)$ as a function of
$x$ or $y$, the resulting function has supremum norm less than or
equal to $\|\mu\|_*$ or $\|\nu\|_*$ times $\|f\|_\infty$, as
appropriate.  For both definitions of $(\mu \times \nu)(f)$ one gets
\begin{equation}
	|(\mu \times \nu)(f)| \le \|\mu\|_* \, \|\nu\|_* \, \|f\|_\infty.
\end{equation}

	If $f(x, y)$ is of the form $f_1(x) \, f_2(y)$, where $f_1$,
$f_2$ are continuous functions on the real line with compact support,
then both definitions of $(\mu \times \nu)(f)$ are equal to $\mu(f_1)
\, \nu(f_2)$.  The two definitions of $\mu \times \nu$ agree on
functions $f$ which are sums of products of this type, and hence on
all continuous functions on the plane with bounded support, by
approximation arguments.  Furthermore,
\begin{equation}
	\sup \{|(\mu \times \nu)(f)| : f \in \mathcal{C}_{00}({\bf R}^2),
					\|f\|_{\infty, {\bf R}^2} \}
		= \|\mu\|_* \, \|\nu\|_*,
\end{equation}
where $\mathcal{C}_{00}({\bf R}^2)$ is the space of continuous
complex-valued functions on the plane with bounded support, and
$\|f\|_{\infty, {\bf R}^2}$ is the supremum norm for such functions.

	Just as for measures on ${\bf R}$, one can consider the notion
of a regular extension for $\mu \times \nu$ to bounded continuous
functions on ${\bf R}^2$.  One can also start with the regular
extensions of $\mu$, $\nu$ to bounded continuous functions on the real
line and use them to deal with $(\mu \times \nu)(f)$ when $f$ is a
bounded continuous function on ${\bf R}^2$.  If $\mu$, $\nu$ have
bounded support in ${\bf R}$, then the product $\mu \times \nu$ has
bounded support in ${\bf R}^2$, and $(\mu \times \nu)(f)$ makes sense
for arbitrary continuous functions $f$ on the plane.  In general there
are approximations of $\mu$, $\nu$ by measures with bounded support,
which lead to approximations of $\mu \times \nu$ by measures with
bounded support.  At any rate, $\mu \times \nu$ has a version of a
regular extension to bounded continuous functions on the plane.

	The \emph{convolution}\index{convolutions} $\mu * \nu$ of
$\mu$, $\nu$ is the measure on ${\bf R}$ defined by saying that $(\mu
* \nu)(f)$ is equal to $\mu \times \nu$ applied to $f(x + y)$ as a
continuous function on the plane.  In particular, this is a bounded
linear functional, and more precisely $\|\mu * \nu\|_* \le \|\mu\|_*
\, \|\nu\|_*$.  The Fourier transform of $\mu * \nu$ is equal to the
product of the Fourier transforms of $\mu$ and $\nu$.

\section{Smooth functions}
\label{smooth functions}
\setcounter{equation}{0}

	Let $\mathcal{C}^\infty({\bf
R})$\index{$Cinfty(R)$@$\mathcal{C}^\infty({\bf R})$} be the space of
complex-valued functions on the real line which are continuously
differentiable of all orders, and let $\mathcal{C}^\infty_{00}({\bf
R})$\index{$Cinfty(R)_00$@$\mathcal{C}^\infty_{00}({\bf R})$} be the
subspace of $\mathcal{C}^\infty({\bf R})$ consisting of functions with
bounded support.  For instance, the function defined by $\exp (-1/x)$
when $x > 0$ and equal to $0$ when $x \le 0$ is in
$\mathcal{C}^\infty({\bf R})$, and one can use this to get nontrivial
functions in $\mathcal{C}^\infty_{00}({\bf R})$.  The \emph{Schwartz
class}\index{Schwartz class} $\mathcal{S}({\bf
R})$\index{$S(R)$@$\mathcal{S}({\bf R})$} consists of the functions $f
\in \mathcal{C}^\infty({\bf R})$ such that $f$ and all of its
derivatives are rapidly decreasing in the sense that for each pair of
nonnegative integers $j$, $l$ there is a $C(j, l) \ge 0$ such that
\begin{equation}
	|f^{(j)}(x)| \le \frac{C(j, l)}{(1 + |x|)^l}
\end{equation}
for every $x \in {\bf R}$, where $f^{(j)}$ denotes the $j$th order
derivative of $f$, with $f^{(0)} = f$.  In particular,
\begin{equation}
	\mathcal{C}^\infty_{00}({\bf R}) \subseteq \mathcal{S}({\bf R}).
\end{equation}

	If $f \in \mathcal{S}({\bf R})$, then $f$ is an integrable
continuous function, and its Fourier transform $\widehat{f}$ is a
bounded continuous function.  One can show that the Fourier transform
of $f$ is continuously differentiable of all orders, and more
precisely that the $l$th derivative of $\widehat{f}$ is equal to
$(-i)^l$ times the Fourier transform of $x^l \, f(x)$ for every
positive integer $l$.  Similarly, $\xi^l \, \widehat{f}(\xi)$ is equal
to $i^l$ times the Fourier transform of the $l$th derivative of $f$
for each positive integer $l$ and hence is bounded.  Continuing in
this way one can show that
\begin{equation}
	\widehat{f} \in \mathcal{S}({\bf R}).
\end{equation}
Similarly, the inverse Fourier transform maps the Schwartz class into
itself, and it follows that the Fourier transform is a one-to-one
mapping of $\mathcal{S}({\bf R})$ onto itself.

	For each positive integer $n$, let $\mathcal{C}^\infty({\bf
T}^n)$\index{$Cinfty(T^n)$@$\mathcal{C}^\infty({\bf T}^n)$} be the
space of complex-valued functions on the $n$-dimensional torus which
are continuously-differentiable of all orders.  One can show that the
Fourier coefficients $a_\alpha$ of such a function decay rapidly in
the sense that for every positive integer $l$ there is a $C(l) \ge 0$
such that
\begin{equation}
	|a_\alpha| \le \frac{C(l)}{(1 + |\alpha_1| + \cdots + |\alpha_n|)^l}
\end{equation}
for every $\alpha \in {\bf Z}^n$, by expressing the Fourier
coefficients of derivatives of $f$ in terms of the $a_\alpha$'s times
products of powers of the $\alpha_j$'s.  Conversely, if $a_\alpha$,
$\alpha \in {\bf Z}^n$, is a family of complex numbers which are
rapidly decreasing in this sense, then $\sum_{\alpha \in {\bf Z}^n}
a_\alpha \, z^\alpha$ is continuously differentiable of all orders on
${\bf T}^n$, and the derivatives are given by differentiating the sum
term-by-term.

	Suppose that $\Phi \in \mathcal{S}({\bf R})$ and $\Phi(0) =
1$.  Put $\phi = \widehat{\Phi}$ and
\begin{equation}
 \phi_\eta(w)
	= \frac{1}{2 \, \pi \, \eta} \, \phi \Big(\frac{w}{\eta} \Big),
\end{equation}
which is $1 / (2 \, \pi)$ times the Fourier transform of $\Phi(\eta \,
x)$ at $w$.  For each $\eta > 0$,
\begin{equation}
	\int_{\bf R} \phi_\eta(w) \, dw = 1,
\end{equation}
because $\Phi(0) = 1$.  For every integrable continuous function $f$
on the real line, consider
\begin{equation}
\label{frac{1}{2 pi} int_R widehat{f}(xi) Phi(eta xi) exp (xi x i) d xi}
 \frac{1}{2 \, \pi} \int_{\bf R} \widehat{f}(\xi) \, \Phi(\eta \, \xi)
					\, \exp (\xi \, x \, i) \, d\xi
\end{equation}
as an extension of the Abel sums\index{Abel summability}
(\ref{frac{1}{2 pi} int_{bf R} widehat{f}(xi) A_eta(xi) exp (xi x i) d
xi}) for the inverse Fourier transform (\ref{frac{1}{2 pi} int_{bf R}
widehat{f}(xi) exp (xi x i) d xi}) applied to $\widehat{f}$.  As
usual, the limit of (\ref{frac{1}{2 pi} int_R widehat{f}(xi) Phi(eta
xi) exp (xi x i) d xi}) as $\eta \to 0$ is equal to (\ref{frac{1}{2
pi} int_{bf R} widehat{f}(xi) exp (xi x i) d xi}) when $\widehat{f}$
is integrable.  In general, we can rewrite (\ref{frac{1}{2 pi} int_R
widehat{f}(xi) Phi(eta xi) exp (xi x i) d xi}) as a double integral
using the definition of $\widehat{f}$.  By interchanging the order of
integration, we get that (\ref{frac{1}{2 pi} int_R widehat{f}(xi)
Phi(eta xi) exp (xi x i) d xi}) is equal to
\begin{equation}
\label{int_R f(y) phi_eta(y - x) dy}
	\int_{\bf R} f(y) \, \phi_\eta(y - x) \, dy.
\end{equation}
As in the previous situation, one can show that the limit of
(\ref{int_R f(y) phi_eta(y - x) dy}) as $\eta \to 0$ is equal to
$f(x)$.

\section{Gaussians}
\label{gaussians}
\setcounter{equation}{0}

	For each positive real number $a$, the corresponding
\emph{Gaussian}\index{Gaussians} function is defined by
\begin{equation}
	G_a(x) = \exp (- a \, x^2).
\end{equation}
This is an integrable continuous function on the real line which is an
element of the Schwartz class $\mathcal{S}({\bf R})$.

	It is a well-known fact from integral calculus that
\begin{equation}
	\int_{\bf R} \exp (- \pi \, x^2) \, dx = 1.
\end{equation}
The trick for showing this is to observe that the integral is a
positive real number whose square can be expressed as the double
integral
\begin{equation}
	\int_{{\bf R}^2} \exp (- \pi \, (x^2 + y^2)) \, dx \, dy.
\end{equation}
Using polar coordinates this double integral reduces to
\begin{equation}
 \int_0^{2 \, \pi} \int_0^\infty \exp (- \pi \, r^2) \, r \, dr \, d\theta
	= \int_0^\infty 2 \, \pi \, r \, \exp (- \pi \, r^2) \, dr = 1,
\end{equation}
since the derivative of $\exp (-\pi \, r^2)$ is equal to $- 2 \, \pi
\, r \, \exp (- \pi \, r^2)$.

	It follows that
\begin{equation}
	\int_{\bf R} \exp (-a \, x^2) \, dx = \sqrt{\frac{\pi}{a}}.
\end{equation}
Specifically, one can use the change of variables $x = \sqrt{(\pi /
a)} \, w$ to get this from the case where $a = \pi$.

	For every $b \in {\bf R}$,
\begin{equation}
	\int_{\bf R} \exp (- a \, (x + (2 \, a)^{-1} \, b)^2) \, dx 
		= \sqrt{\frac{\pi}{a}},
\end{equation}
by translation-invariance of the integral.  Hence
\begin{equation}
	\int_{\bf R} \exp (- a \, x^2 - b \, x) \, dx
		= \sqrt{\frac{\pi}{a}} \, \exp ((4 \, a)^{-1} \, b^2).
\end{equation}

	This suggests that
\begin{equation}
	\int_{\bf R} \exp (- a \, x^2 - \pi \, \xi \, x \, i) \, dx
		= \sqrt{\frac{\pi}{a}} \, \exp (- (4 \, a)^{-1} \, \xi^2),
\end{equation}
which is to say that
\begin{equation}
	\widehat{G_a}(\xi) = \sqrt{\frac{\pi}{a}} \, G_{(4 \, a)^{-1}}(\xi),
\end{equation}
$\xi \in {\bf R}$.

	One can show this rigorously using complex analysis, through
Cauchy's theorem or analytic continuation in $b$.

	Alternatively one can use differential equations, observing a
linear relation between the derivative of $G_a(x)$ and $x \, G_a(x)$
and a similar relationship for the Fourier transform.

\section{Plancherel's theorem}
\label{plancherel's theorem}
\setcounter{equation}{0}

	If $\phi_1$, $\phi_2$ are continuous integrable functions on
the real line with integrable Fourier transforms, then
\begin{equation}
 \int_{\bf R} \widehat{\phi_1}(\xi) \, \overline{\widehat{\phi_2}(\xi)} \, d\xi
   = \int_{\bf R} \int_{\bf R} \widehat{\phi_1}(\xi) \, \overline{\phi_2(x)}
			\, \exp (\xi \, x \, i) \, dx \, d\xi.
\end{equation}
Interchanging the order of integration and using the formula for the
inverse Fourier transform, we get that
\begin{equation}
\label{fourier transform and inner products}
 \int_{\bf R} \widehat{\phi_1}(\xi) \, \overline{\widehat{\phi_2}(\xi)} \, d\xi
	= 2 \, \pi \int_{\bf R} \phi_1(x) \, \overline{\phi_2(x)} \, dx.
\end{equation}
In particular,
\begin{equation}
\label{fourier transforms and square norm}
	\int_{\bf R} |\widehat{\phi}(\xi)|^2 \, d\xi
		= 2 \, \pi \int_{\bf R} |\phi(x)|^2 \, dx
\end{equation}
when $\phi$ is a continuous integrable function on the real line with
integrable Fourier transform.

	More generally, a continuous integrable function on the real
line whose square is integrable has square-integrable Fourier
transform, and the identities (\ref{fourier transform and inner
products}) and (\ref{fourier transforms and square norm}) carry over
to these functions.  Observe that a continuous integrable function
with integrable Fourier transform is bounded and hence square
integrable.

	In order to show this extension, one can regularize the
integrals in the preceding computations as in the Abel
summability\index{Abel summability} techniques employed several times
now.  One can also approximate an integrable and square-integrable
function on the real line simultaneously in the norms $\|f\|_1$ and
$\|f\|_2$ by integrable continuous functions with integrable Fourier
transforms, e.g., by smooth functions with bounded support.

\section{Bounded functions}
\label{bounded functions}
\setcounter{equation}{0}

	If $f$, $\phi$ are continuous integrable functions on the real
line, then
\begin{equation}
	\int_{\bf R} \widehat{f}(\xi) \, \phi(\xi) \, d\xi
		= \int_{\bf R} f(x) \, \widehat{\phi}(x) \, dx,
\end{equation}
basically because they are both equal to the double integral
\begin{equation}
	\int_{\bf R} \int_{\bf R} f(x) \, \phi(\xi) \, \exp (- \xi \, x \, i).
\end{equation}
Let $\mathcal{E}({\bf R})$\index{$E(R)$@$\mathcal{E}({\bf R})$} be the
vector space of continuous integrable functions $\phi$ on ${\bf R}$
such that the Fourier transform $\widehat{\phi}$ of $\phi$ is
integrable too, which implies that $\phi, \widehat{\phi} \in
\mathcal{C}_0({\bf R})$.  If $f$ is a bounded continuous function on
the real line, then
\begin{equation}
	L_f(\phi) = \int_{\bf R} f(x) \, \widehat{\phi}(x) \, dx
\end{equation}
defines a linear mapping from $\mathcal{E}({\bf R})$ into the complex
numbers which is the same as
\begin{equation}
	L_f(\phi) = \int_{\bf R} \widehat{f}(\xi) \, \phi(\xi) \, d\xi
\end{equation}
when $f$ is integrable.

	We can think of $L_f$ as a kind of \emph{generalized Fourier
transform}\index{generalized Fourier transforms} which makes sense for
bounded continuous functions on the real line.  For example, if
\begin{equation}
	f(x) = \exp (a \, x \, i)
\end{equation}
for some $a \in {\bf R}$, then
\begin{equation}
	L_f(\phi) = 2 \, \pi \, \phi(a)
\end{equation}
for every $\phi \in \mathcal{E}({\bf R})$.  As a version of Abel
sums\index{Abel summability} for the inverse Fourier transform of a
function $f \in \mathcal{BC}({\bf R})$, we can apply $L_f$ to
\begin{equation}
	\frac{1}{2 \, \pi} \, A_\eta(\xi) \, \exp (\xi \, v \, i)
\end{equation}
as a function of $\xi$ for every $v \in {\bf R}$ and $\eta > 0$.  As
usual this is equal to the Poisson integral $P_\eta(f)(v)$ of $f$, and
we recover $f(v)$ as $\eta \to 0$.  One can also use other functions
instead of $A_\eta(\xi)$ as in Section \ref{smooth functions}.

	Let $\sigma(x)$ be the function on the real line which is
equal to $-1$ when $x < 0$, to $0$ when $x = 0$, and to $+1$ when $x >
0$.  Although this bounded function is not quite continuous,
\begin{equation}
	L_\sigma(\phi) = \int_{\bf R} \sigma(x) \, \widehat{\phi}(x) \, dx
\end{equation}
still makes sense as a generalized Fourier transform of $\sigma$, and
it will be convenient for us to restrict our attention now to $\phi
\in \mathcal{S}({\bf R})$.  For each $\eta > 0$, put
\begin{equation}
	L_{\sigma, \eta}(\phi) 
	 = \int_{\bf R} \sigma(x) \, A_\eta(x) \, \widehat{\phi}(x) \, dx,
\end{equation}
where $A_\eta(x) = \exp ( - \eta \, |x|)$.  Clearly
\begin{equation}
	\lim_{\eta \to 0} L_{\sigma, \eta}(\phi) = L_\sigma(\phi)
\end{equation}
for every $\phi \in \mathcal{S}({\bf R})$.

	Because
\begin{equation}
	B_\eta(x) = \sigma(x) \, A_\eta(x)
\end{equation}
is integrable, with a jump discontinuity at $x = 0$, its Fourier
transform can be defined in the usual way, and
\begin{equation}
	L_{\sigma, \eta}(\phi)
		= \int_{\bf R} \widehat{B_\eta}(\xi) \, \phi(\xi) \, d\xi.
\end{equation}
Observe that
\begin{eqnarray}
	\quad \widehat{B_\eta}(\xi)
		& = & -\int_{-\infty}^0 \exp ((\eta - \xi \, i) \, x) \, dx
			+ \int_0^\infty \exp (-(\eta + \xi \, i) \, x) \, dx
								\\
		& = & \frac{-1}{\eta - \xi \, i} + \frac{1}{\eta + \xi \, i}
			= \frac{2 \, \xi \, i}{\eta^2 + \xi^2},
						\nonumber
\end{eqnarray}
which is an odd function,
\begin{equation}
	\widehat{B_\eta}(-\xi) = - \widehat{B_\eta}(\xi).
\end{equation}
We can rewrite $L_{\sigma, \eta}(\phi)$ as
\begin{equation}
	L_{\sigma, \eta}(\phi)
 = \int_{|\xi| \le 1} \widehat{B_\eta}(\xi) \, (\phi(\xi) - \phi(0)) \, d\xi
  + \int_{|\xi| > 1} \widehat{B_\eta}(\xi) \, \phi(\xi) \, d\xi.
\end{equation}
Because of the smoothness and integrability of $\phi$,
\begin{equation}
	\lim_{\eta \to 0} L_{\sigma, \eta}(\phi)
 = 2 \, i \, \int_{|\xi| \le 1} \frac{\phi(\xi) - \phi(0)}{\xi} \, d\xi
  + 2 \, i \, \int_{|\xi| > 1} \frac{\phi(\xi)}{\xi} \, d\xi.
\end{equation}
Equivalently,
\begin{equation}
	\lim_{\eta \to 0} L_{\sigma, \eta}(\phi)
 = \lim_{\epsilon \to 0} 2 \, i \, \int_{|\xi| > \epsilon} 
					\frac{\phi(\xi)}{\xi} \, d\xi.
\end{equation}

\section{Subharmonic functions}
\label{subharmonic functions}
\setcounter{equation}{0}

	Let $u(z)$ be a twice-continuously differentiable real-valued
function defined on an open set in the complex plane.  We say that $u$
is \emph{subharmonic}\index{subharmonic functions} if
\begin{equation}
	\Delta u(z) \ge 0
\end{equation}
for every $z$ in the domain of $u$, where
\begin{equation}
	\Delta
 = \frac{\partial^2}{\partial x^2} + \frac{\partial^2}{\partial y^2}
\end{equation}
is the usual Laplace operator, and $x$, $y$ are the real and imaginary
parts of $z$.  We say that $u$ is \emph{strictly subharmonic} if
\begin{equation}
	\Delta u(z) > 0
\end{equation}
for every $z$ in the domain of $u$.

	For the sake of concreteness let us suppose that $u$ is a
continuous real valued function on the closed unit disk $\{z \in {\bf
C} : |z| \le 1\}$ which is twice-continuously differentiable on the
open unit disk $\{z \in {\bf C} : |z| < 1\}$.  If $u$ is strictly
subharmonic on the open unit disk, then the second-derivative test
from calculus implies that $u$ does not have any local maxima in the
open unit disk.  However, the maximum of $u$ on the closed unit disk
is attained, because the closed unit disk is compact and $u$ is
continuous.  Hence the maximum is attained on the unit circle.

	Therefore
\begin{equation}
	\max \{u(z) : z \in {\bf C}, \ |z| \le 1\}
		= \max \{u(z) : z \in {\bf C}, |z| = 1 \}.
\end{equation}
This also works when $u$ is subharmonic on the open unit disk.  For in
this event $u_\epsilon(z) u(z) + \epsilon \, |z|^2$ is strictly
subharmonic for every $\epsilon > 0$.  By applying the maximum
principle to $u_\epsilon$ and sending $\epsilon \to 0$ we recover the
maximum principle for $u$.

	If $u$ is harmonic on the open unit disk, then $u$ and $-u$
are subharmonic.  In particular, if $u(z) = 0$ for every $z$ in the
unit circle, then $u$ vanishes everywhere.  Equivalently, two
continuous real-valued functions on the closed unit disk which are
harmonic on the open unit disk and are equal at every point in the
unit circle are equal everywhere.  The same statement also holds for
complex-valued functions, by considering the real and imaginary parts.

\appendix

\section{Metric spaces}
\label{metric spaces}
\renewcommand{\thetheorem}{A.\arabic{equation}}
\renewcommand{\theequation}{A.\arabic{equation}}
\setcounter{equation}{0}

	A \emph{metric space}\index{metric spaces} is a nonempty set
$M$ together with a distance function $d(x, y)$ defined for $x, y \in
M$ such that $d(x, y)$ is a nonnegative real number for every $x, y
\in M$,
\begin{equation}
	d(x, y) = 0
\end{equation}
if and only if $x = y$,
\begin{equation}
	d(y, x) = d(x, y)
\end{equation}
for every $x, y \in M$, and
\begin{equation}
	d(x, z) \le d(x, y) + d(y, z)
\end{equation}
for every $x, y, z \in M$.  The last property is known as the triangle
inequality.

	For example, the real line ${\bf R}$ is a metric space with
the standard metric $|x - y|$, and the complex numbers ${\bf C}$ form
a metric space with the standard metric $|z - w|$.  If $V$ is a vector
space equipped with a norm $N(v)$, then $N(v - w)$ defines a metric on
$V$.  On ${\bf C}^n$, $\|v - w\|_p^p$ is a metric when $0 < p \le 1$.
Here $\|v\|_p$ is as in Section \ref{norms}, and is a norm on ${\bf
C}^n$ when $p \ge 1$.  If $(M, d(x, y))$ is any metric space and $E
\subseteq M$, $E \ne \emptyset$, then $E$ is a metric space using the
restriction of $d(x, y)$ to $E$ as the metric.

	Let $(M, d(x, y))$ be a metric space.  For every $x \in M$ and
$r > 0$, the open ball with center $x$ and radius $r$ is given by
\begin{equation}
	B(x, r) = \{y \in M : d(x, y) < r \},
\end{equation}
and the closed ball with center $x$ and radius $r$ is given by
\begin{equation}
	\overline{B}(x, r) = \{y \in M : d(x, y) \le r \}.
\end{equation}

	A set $U \subseteq M$ is said to be \emph{open}\index{open
sets} if for every $x \in U$ there is an $r > 0$ such that $B(x, r)
\subseteq U$.  One can check that open balls are open sets.

	If $A \subseteq M$ and $p \in M$, then $p$ is a \emph{limit
point}\index{limit points} of $A$ if for every $r > 0$ there is a $q
\in A$ such that $q \ne p$ and $d(p, q) < r$.  Similarly, $p$ is an
\emph{accumulation point}\index{accumulation points} of $A$ if for
every $r > 0$ there is a $q \in A$ such that $d(p, q) < r$.  Every
accumulation point of $A$ is an element of $A$, or a limit point of
$A$, or both.

	We say that $E \subseteq M$ is \emph{closed}\index{closed
sets} if every $p \in M$ which is a limit point of $E$ is also an
element of $E$, which is the same as saying that every accumulation
point of $E$ is an element of $E$.  One can check that closed balls
are closed sets.

	The \emph{complement} of a set $A \subseteq M$ is the set $M
\backslash A$ of $x \in M$ such that $x$ is not an element of $A$.
One can show that $A \subseteq M$ is an open set if and only if $M
\backslash A$ is closed.

	The union of any collection of open subsets of $M$ is an open
set.  The intersection of any collection of closed subsets of $M$ is a
closed set.  These statements are easy to check, just from the
definitions, and they also correspond to each other by taking
complements as in the previous paragraph.

	The intersection of finitely many open subsets of $M$ is an
open set.  It follows from the fact about complements that the union
of finitely many closed subsets of $M$ is closed.

	The \emph{closure}\index{closure of a set} of $E \subseteq M$
is denoted $\overline{E}$ and defined to be the set of accumulation
points of $E$ in $M$, which is the same as the union of $E$ and the
set of limit points of $E$ in $M$.  Observe that $E$ is closed if and
only if $\overline{E} = E$.  One can check that $\overline{E}$ is
automatically closed, which is to say that an accumulation point of
the set of accumulation points of $E$ is an accumulation point of $E$.

	We say that $E \subseteq M$ is \emph{bounded}\index{bounded sets}
if there are $p \in M$ and $r > 0$ such that
\begin{equation}
	E \subseteq B(p, r),
\end{equation}
in which event the \emph{diameter} of $E$ is given by
\begin{equation}
	\diam E = \sup \{d(x, y) : x, y \in E \},
\end{equation}
with $\diam E = 0$ when $E = \emptyset$.  The closure of a bounded set
is bounded and has the same diameter.

	Let $\{x_j\}_{j=1}^\infty$ be a sequence of points in $M$.  We
say that $\{x_j\}_{j=1}^\infty$ \emph{converges}\index{convergent
sequences} to $x \in M$ if for each $\epsilon > 0$ there is an $L \ge
1$ such that
\begin{equation}
	d(x_j, x) < \epsilon
\end{equation}
for every $j \ge L$.  We say that $\{x_j\}_{j=1}^\infty$ is a
\emph{Cauchy sequence}\index{Cauchy sequences} if for each $\epsilon >
0$ there is an $L \ge 1$ such that
\begin{equation}
	d(x_j, x_l) < \epsilon
\end{equation}
for every $j, l \ge L$.  It is easy to see that convergent sequences
are Cauchy sequences.  If $\{x_j\}_{j=1}^\infty$ converges to $x$ in
$M$, then $x$ is said to be the \emph{limit} of the sequence, also
expressed by
\begin{equation}
	\lim_{j \to \infty} x_j = x,
\end{equation}
and one can check that the limit $x$ is unique.

	A point $p \in M$ is an accumulation point of $E \subseteq M$
if and only if there is a sequence of elements of $E$ which converges
to $p$.  Hence $E \subseteq M$ is closed if and only if every sequence
of elements of $E$ which converges in $M$ has its limit in $E$.

	If $\{z_j\}_{j=1}^\infty$, $\{w_j\}_{j=1}^\infty$ are
sequences of real or complex numbers which converge to the real or
complex numbers $z$, $w$, respectively, then the sequences
\begin{equation}
	\{z_j + w_j\}_{j=1}^\infty, \quad \{z_j \, w_j\}_{j=1}^\infty
\end{equation}
of sums and products converge to $z + w$, $z \, w$, respectively.

	More generally, if $V$ is a complex vector space equipped with
a norm and thus a metric, and if $\{v_j\}_{j=1}^\infty$,
$\{w_j\}_{j=1}^\infty$ are sequences of vectors in $V$ converging to
$v, w \in V$, then
\begin{equation}
	\lim_{j \to \infty} v_j + w_j = v + w.
\end{equation}
If $\{\alpha_j\}_{j=1}^\infty$ is a sequence of complex numbers
converging to $\alpha \in {\bf C}$ and $\{v_j\}_{j=1}^\infty$ is a
sequence of vectors in $V$ converging to $v \in V$, then
\begin{equation}
	\lim_{j \to \infty} \alpha_j \, v_j = \alpha \, v.
\end{equation}

	A metric space is said to be
\emph{complete}\index{completeness} if every Cauchy sequence in the
space converges.  The real and complex numbers are complete with
respect to their standard metrics.

	If $(M, d(x, y))$ is a metric space, then $E \subseteq M$ is
said to be \emph{dense} in $M$ if $\overline{E} = M$.  Equivalently,
$E$ is dense in $M$ if for every $x \in M$ there is a sequence
$\{x_j\}_{j=1}^\infty$ of elements of $E$ which converges to $x$.  For
example, the rational numbers are dense in ${\bf R}$.

\section{Compact sets}
\label{Compact sets}
\renewcommand{\thetheorem}{B.\arabic{equation}}
\renewcommand{\theequation}{B.\arabic{equation}}
\setcounter{equation}{0}

	Let $\{x_j\}_{j=1}^\infty$ be a sequence with terms in any
set.  A \emph{subsequence}\index{subsequences} of
$\{x_j\}_{j=1}^\infty$ is a sequence of the form
$\{x_{j_l}\}_{l=1}^\infty$, where $\{j_l\}_{l=1}^\infty$ is a strictly
increasing sequence of positive integers.  In particular,
$\{x_j\}_{j=1}^\infty$ is a subsequence of itself.

	In a metric space $(M, d(x, y))$, every subsequence of a
convergent sequence converges to the same point in $M$.  If a Cauchy
sequence in $M$ has a convergent subsequence, then the Cauchy sequence
converges to the same point in $M$.

	A set $E \subseteq M$ is said to be \emph{sequentially
compact}\index{sequential compactness} if every sequence of elements
of $E$ has a subsequence which converges to an element of $E$.  In
particular, if a sequence of elements of $E$ converges in $M$, then
the limit has to be in $E$.  Hence a sequentially compact set is
closed.

	As a partial converse, if $E \subseteq M$ is closed, $E_1
\subseteq M$ is sequentially compact, and $E \subseteq E_1$, then $E$
is sequentially compact.  For sequential compactness of $E_1$ implies
that every sequence in $E$ has a convergent subsequence, whose limit
is in $E$ since $E$ is closed.

	Sequentially compact sets are bounded.  For if $E \subseteq M$
is not bounded and $p \in M$, then for each positive integer $j$ there
is an $x_j \in E$ such that $d(x_j, p) \ge j$, and it is easy to see
that $\{x_j\}_{j=1}^\infty$ has not convergent subsequence.

	A set $E \subseteq M$ is said to be \emph{totally
bounded}\index{totally bounded sets} if for each $\epsilon > 0$ there
are finitely many points $p_1, \ldots, p_n \in E$ such that $E$ is
contained in the union of the balls $B(p_1, \epsilon), \ldots, B(p_n,
\epsilon)$.  If $E$ is not totally bounded, then there is a sequence
$\{x_j\}_{j=1}^\infty$ of elements of $E$ such that $d(x_j, x_l) \ge
\epsilon$ when $j \ne l$.  It follows that sequentially compact sets
are totally bounded.

	If a sequence has its terms in the union of two sets, then
there is a subsequence whose terms are all in one of the sets.  It
follows that the union of two sequentially compact sets is
sequentially compact.

	A set $E \subseteq M$ has the \emph{limit point
property}\index{limit point property} if every infinite set $A
\subseteq E$ has a limit point contained in $E$.  If $A \subseteq E$
is infinite, then there is a sequence $\{a_j\}_{j=1}^\infty$ with $a_j
\in A$ for all $j$ and $a_j \ne a_l$ when $j \ne l$, and the limit of
any convergent subsequence of this sequence is a limit point of $A$.
Hence sequential compactness implies the limit point property.
Conversely, if $\{x_j\}_{j=1}^\infty$ is a sequence of elements of $E$
and $A$ is the set of all the $x_j$'s, then either $A$ has only
finitely many elements and the sequence has a constant subsequence, or
$A$ is infinite and one can check that any limit point of $A$ is the
limit of a subsequence of $\{x_j\}_{j=1}^\infty$.  Therefore the limit
point propert implies sequential compactness.

	In these notes we shall use the term ``compact''\index{compact
sets} to refer to a set which is sequentially compact or has the limit
point property.  Equivalently, $E \subseteq M$ is compact if every
covering of $E$ by open subsets of $M$ can be reduced to a covering by
finitely many of the open sets, but we shall not use this here.

	Suppose that $E \subseteq M$ has the property that every
sequence of elements of $E$ has a subsequence which is a Cauchy
sequence.  The same argument as above shows that $E$ is totally
bounded.  Conversely, one can show that a totally bounded set has this
property.

	As a consequence, if $(M, d(x, y))$ is a complete metric
space, and $E \subseteq M$ is closed and totally bounded, then $E$ is
sequentially compact.  The converse works in any metric space by the
earlier remarks.

\section{Continuous functions}
\label{continuous functions}
\renewcommand{\thetheorem}{C.\arabic{equation}}
\renewcommand{\theequation}{C.\arabic{equation}}
\setcounter{equation}{0}

	Let $(M, d(x, y))$ be a metric space, and let $f$ be a
complex-valued function on $M$.  We say that $f$ is
\emph{continuous}\index{continuous functions} at a point $x \in M$ if
for every $\epsilon > 0$ and $x \in M$ there is a $\delta > 0$ such
that
\begin{equation}
\label{|f(x) - f(y)| < epsilon}
	|f(x) - f(y)| < \epsilon
\end{equation}
for each $y \in M$ such that $d(x, y) < \delta$.  Equivalently, $f$ is
continuous at $x$ if for every sequence $\{x_j\}_{j=1}^\infty$ of
elements of $M$ which converges to $x$ we have that
\begin{equation}
	\lim_{j \to \infty} f(x_j) = f(x).
\end{equation}

	The space of complex-valued functions on $M$ which are
continuous at every point in $M$ is denoted $\mathcal{C}(M)$.
Constant functions are obviously continuous, and one can show that
sums and products of continuous functions are continuous, which
implies that $\mathcal{C}(M)$ is a commutative algebra over ${\bf C}$.

	A continuous complex-valued function $f$ on $M$ is said to be
\emph{bounded}\index{bounded functions} if there is an $A \ge 0$ such
that
\begin{equation}
	|f(x)| \le A
\end{equation}
for every $x \in M$.  The space of bounded continuous complex-valued
functions on $M$ is denoted $\mathcal{BC}(M)$ and is a subalgebra of
$\mathcal{C}(M)$.  The supremum norm of $f \in \mathcal{BC}(M)$ is given by
\begin{equation}
	\|f\|_\infty = \|f\|_{\infty, M} = \sup \{|f(x)| : x \in M \},
\end{equation}
which one can check is a norm on $\mathcal{BC}(M)$ as a complex vector
space and also satisfies
\begin{equation}
	\|f_1 \, f_2\|_\infty \le \|f_1\|_\infty \, \|f_2\|_\infty
\end{equation}
for every $f_1, f_2 \in \mathcal{BC}(M)$.

	If $f$ is a continuous complex-valued function on $M$, $E
\subseteq M$ is sequentially compact, and $\{x_j\}_{j=1}^\infty$ is
any sequence of elements of $E$, then there is a subsequence
$\{x_{j_l}\}_{l=1}^\infty$ of $\{x_j\}_{j=1}^\infty$ which converges
to a point $x \in E$, and $\{f(x_j)\}_{j=1}^\infty$ converges in ${\bf
C}$ to $f(x)$.  One can use this to show that $f$ is bounded on $E$
and $|f|$ attains its maximum on $E$ when $E \ne \emptyset$.

	A complex-valued function $f$ on $M$ is said to be
\emph{uniformly continuous}\index{uniform continuity} if for each
$\epsilon > 0$ there is a $\delta > 0$ such that (\ref{|f(x) - f(y)| <
epsilon}) holds for every $x, y \in M$ with $d(x, y) < \delta$.
Uniformly continuous functions are automatically continuous.  The
space of uniformly continuous complex-valued functions on $M$ is
denoted $\mathcal{UC}(M)$ and is a linear subspace of
$\mathcal{C}(M)$.

	Suppose that $f$ is a continuous function on $M$ which is not
uniformly continuous.  Then there is an $\epsilon > 0$ and sequences
$\{x_j\}_{j=1}^\infty$, $\{y_j\}_{j=1}^\infty$ of elements of $M$ such
that
\begin{equation}
	\lim_{j \to \infty} d(x_j, y_j) = 0
\end{equation}
and
\begin{equation}
	|f(x_j) - f(y_j)| \ge \epsilon.
\end{equation}

	If $M$ is sequentially compact, then every continuous function
on $M$ is uniformly continuous.  Otherwise there would be a
subsequence $\{x_{j_l}\}_{l=1}^\infty$ of $\{x_j\}_{j=1}^\infty$ as in
the previous paragraph which converges to a point $x \in M$, the
corresponding subsequence $\{y_j\}_{j=1}^\infty$ would also converge
to $x$, and continuity of $f$ at $x$ would imply that
$\{f(x_{j_l})\}_{l=1}^\infty$, $\{f(y_{j_l})\}_{l=1}^\infty$ converge
to $f(x)$ as sequences of complex numbers, a contradiction.

	Let $\{f_j\}_{j=1}^\infty$ be a sequence of complex-valued
continuous functions on $M$ which converges
\emph{uniformly}\index{uniform convergence} to a complex-valued
function $f$ on $M$ in the sense that for each $\epsilon > 0$ there is
an $L \ge 1$ such that
\begin{equation}
	|f_j(x) - f(x)| < \epsilon
\end{equation}
for every $x \in M$ and $j \ge L$.  In this event one can show that
the limiting function $f$ is also a continuous function on $M$.
Similarly, if the $f_j$'s are uniformly continuous, then $f$ is
uniformly continuous.

	If the $f_j$'s are bounded, then $f$ is bounded too.  For
bounded continuous functions, uniform convergence is equivalent to
convergence in $\mathcal{BC}(M)$ with respect to the supremum metric.

	One can check that $\mathcal{BC}(M)$ is complete with respect
to the supremum metric.  For if $\{f_j\}_{j=1}^\infty$ is a Cauchy
sequence in $\mathcal{BC}(M)$ with respect to the supremum metric,
then $\{f_j(x)\}_{j=1}^\infty$ is a Cauchy sequence of complex numbers
for every $x \in M$.  The completeness of the complex numbers implies
that $\{f_j(x)\}_{j=1}^\infty$ converges in ${\bf C}$ for every $x \in
M$, and $f(x)$ denotes the limit, then one can check that
$\{f_j\}_{j=1}^\infty$ converges uniformly to $f$ on $M$.

	Let $\mathcal{BUC}(M)$ be the space of bounded uniformly
continuous complex-valued functions on $M$, i.e.,
\begin{equation}
	\mathcal{BUC}(M) = \mathcal{BC}(M) \cap \mathcal{UC}(M).
\end{equation}
One can check that $\mathcal{BUC}(M)$ is a closed subalgebra of
$\mathcal{BC}(M)$.

	If $f$ is a continuous complex-valued function on $M$, then
the zero set of $f$, consisting of $x \in M$ such that $f(x) = 0$, is
a closed set.  If $f_1$, $f_2$ are continuous complex-valued functions
on $M$, then the set of $x \in M$ such that $f_1(x) = f_2(x)$ is
closed, since this is the same as the zero set of $f_1 - f_2$.  In
particular, if two continuous functions on $M$ are equal at every
element of a dense set, then they are equal at every point in $M$.

	If $E \subseteq M$ is dense and $f$ is a uniformly continuous
complex-valued function on $E$, then there is an extension of $f$ to a
uniformly continuous complex-valued function on $M$.  The main point
is that for each $x \in M$ there is a sequence $\{x_j\}_{j=1}^\infty$
of elements of $E$ which converges to $x$ in $M$, and which is a
Cauchy sequence as a sequence in $E$.  Uniform continuity of $f$ on
$E$ implies that $\{f(x_j)\}_{j=1}^\infty$ is a Cauchy sequence and
hence converges as a sequence of complex numbers.  If
$\{y_j\}_{j=1}^\infty$ is another sequence of elements of $E$ which
converges to $x$, then uniform continuity of $f$ on $E$ also implies
that
\begin{equation}
	\lim_{j \to \infty} f(x_j) = \lim_{j \to \infty} f(y_j).
\end{equation}
For every $x \in M$, the value of the extension of $f$ at $x$ is
defined to be any such limit, and it is easy to check that uniform
continuity carries over to the extension.

\section{Lipschitz functions}
\label{lipschitz functions}
\renewcommand{\thetheorem}{D.\arabic{equation}}
\renewcommand{\theequation}{D.\arabic{equation}}
\setcounter{equation}{0}

	Let $(M, d(x, y))$ be a metric space.  A complex-valued
function $f$ on $M$ is said to be \emph{Lipschitz}\index{Lipschitz
functions} if there is a $C \ge 0$ such that
\begin{equation}
	|f(x) - f(x)| \le C \, d(x, y)
\end{equation}
for every $x, y \in M$.  If $f$ is real-valued, then this condition is
equivalent to
\begin{equation}
	f(x) \le f(y) + C \, d(x, y)
\end{equation}
for every $x, y \in M$.  In particular, this holds for $f_p(x) = d(x,
p)$ for every $p \in M$ with $C = 1$, by the triangle inequality.

	More generally, if $a$ is a positive real number, then a
complex-valued function $f$ on $M$ is said to be \emph{Lipschitz of
order $a$} if there is a $C \ge 0$ such that
\begin{equation}
	|f(x) - f(y)| \le C \, d(x, y)^a.
\end{equation}
A Lipschitz function of any order $a > 0$ is uniformly continuous.

	On any metric space a constant function is Lipschitz of order
$a$ for each $a > 0$.  On the real line or the $n$-dimensional torus,
for instance, one can show that a Lipschitz function of order $a > 1$
is constant.

	The sum of two Lipschitz functions of order $a$ is a Lipschitz
function of order $a$, and the product of a Lipschitz function of
order $a$ with a constant is a Lipschitz function of order $a$.  The
product of two bounded Lipschitz functions of order $a$ is a Lipschitz
function of order $a$.  If $f_1$, $f_2$ are two real-valued Lipschitz
functions of order $a$ on $M$, both with constant $C$, then $\max
(f_1, f_2)$ and $\min (f_1, f_2)$ are Lipschitz functions on $M$ of
order $a$ and with constant $C$.

	When $0 < a \le 1$,
\begin{equation}
	(r + t)^a \le r^a + t^a
\end{equation}
for every $r, t \ge 0$, because
\begin{eqnarray}
	r + t & \le & \max(r, t)^{1 - a} \, (r^a + t^a) 	\\
	 	& \le & (r^a + t^a)^{(1 - a) / a} \, (r^a + t^a)
						\nonumber \\
		& = & (r^a + t^a)^{1/a}.	\nonumber
\end{eqnarray}
This implies that $d(x, y)^a$ satisfies the triangle inequality and is
therefore a metric on $M$, and a Lipschitz function on $M$ of order
$a$ with respect to $d(x, y)$ is the same as a Lipschitz function on
$M$ of order $1$ with respect to $d(x, y)^a$.

	Let us restrict our attention now to Lipschitz functions of
order $1$.

	Suppose that $f$ is a bounded continuous real-valued function
on $M$.  For each positive integer $j$, put
\begin{equation}
	f_j(x) = \inf \{f(y) + j \, d(x, y) : y \in M \}.
\end{equation}

	If $c \in {\bf R}$ and $f(w) \ge c$ for every $w \in M$, then
$f(y) + j \, d(x, y) \ge c$ for every $x, y \in M$ and $j \ge 1$, and
hence
\begin{equation}
	f_j(x) \ge c
\end{equation}
for every $x \in M$ and $j \ge 1$.  Since we can take $y = x$ in the
infimum in the definition of $f_j(x)$, we get that
\begin{equation}
	f_j(x) \le f(x)
\end{equation}
for every $x \in M$ and $j \ge 1$.

	For each $j$, $f_j(x)$ is a Lipschitz function on $M$ with
constant $j$.  One can check this using the fact that $f(y) + j \,
d(x, y)$ is Lipschitz with constant $j$ as a function of $x$ for every
$y \in M$ and positive integer $j$.

	Because $f$ is bounded, only $y \in M$ with $d(x, y) = O(1/j)$
are important in the definition of $f_j(x)$, and one can use this and
the continuity of $f$ to show that
\begin{equation}
	\lim_{j \to \infty} f_j(x) = f(x)
\end{equation}
for every $x \in M$.  If $f$ is uniformly continuous, then one can
show that the convergence is uniform.  By applying this to the real
and imaginary parts of bounded uniformly continuous complex-valued
functions on $M$, it follows that the bounded Lipschitz functions are
dense in $\mathcal{BUC}(M)$.

\section{Ultrametric spaces}
\label{ultrametric spaces}
\renewcommand{\thetheorem}{E.\arabic{equation}}
\renewcommand{\theequation}{E.\arabic{equation}}
\setcounter{equation}{0}

	An \emph{ultrametric space}\index{ultrametric spaces} is a
metric space $(M, d(x, y))$ in which the distance function $d(x, y)$
satisfies the stronger version of the triangle inequality,
\begin{equation}
	d(x, z) \le \max (d(x, y), d(y, z))
\end{equation}
for every $x, y, z \in M$.  In this case, $d(x, y)^a$ is an
ultrametric on $M$ for every $a > 0$, and these ultrametrics determine
the same topology on $M$.

	Let $(M, d(x, y))$ be an ultrametric space, and consider the
open ball $B(x, r)$ for some $x \in M$ and $r > 0$.  Because of the
ultrametric version of the triangle inequality, if $d(x, y) < r$ and
$d(x, z) \ge r$, then $d(y, z) \ge r$, which implies that $B(x, r)$ is
a closed set.  Similarly, for every $y \in M$ which satisfies $d(x, y)
\le r$, $\overline{B}(y, r)$ is contained in $\overline{B}(x, r)$,
and therefore $\overline{B}(x, r)$ is an open set in $M$.

	On any nonempty set the discrete metric, which assigns
distance $1$ to every pair of distinct points, is an ultrametric.  As
a more complicated class of examples, let $A$ be a nonempty set, and
let $\Sigma(A)$ be the set of sequences $\{x_j\}_{j=1}^\infty$ with
$x_j \in A$ for each $j$.  For $0 < \rho \le 1$, put $d_\rho(x, y)$
equal to $0$ when $x = y$, and put
\begin{equation}
	d_\rho(x, y) = \rho^l
\end{equation}
when $x_j = y_j$ for $j \le l$ and $x_{l+1} \ne y_{l + 1}$, $x, y \in
\Sigma(A)$, $l \ge 0$.  One can check that $d_\rho(x, y)$ defines an
ultrametric on $\Sigma(A)$ which is equal to the discrete metric when
$\rho = 1$, and which satisfy
\begin{equation}
	d_{\rho^a}(x, y) = d_\rho(x, y)^a
\end{equation}
for each $a > 0$.  One can also show that $\Sigma(A)$ is compact when
$0 < \rho < 1$ if and only if $A$ has only finitely many elements.

	Let $\mathcal{B}$ be the space of all binary sequences, which
is the same as $\Sigma(A)$ with $A = \{0, 1\}$.  There is a mapping
from $\mathcal{B}$ onto the unit interval $[0, 1]$ in the real line
defined by
\begin{equation}
	x = \{x_j\}_{j=1}^\infty \mapsto \sum_{j=1}^\infty x_j \, 2^{-j}.
\end{equation}
This mapping is Lipschitz with respect to the metric $d_\rho(x, y)$
described in the previous paragraph with $\rho = 1/2$.  There is a
one-to-one mapping from $\mathcal{B}$ onto the classical Cantor
middle-thirds set in the real line defined by
\begin{equation}
	x = \{x_j\}_{j=1}^\infty \mapsto \sum_{j=1}^\infty 2 \, x_j \, 3^{-j}.
\end{equation}
This mapping is Lipschitz and moreover bi-Lipschitz in the sense that
distances in the domain and the corresponding distances in the image
are each bounded by a constant multiple of the other when we use the
metric described in the previous paragraph with $\rho = 1/3$.

\newpage

\addcontentsline{toc}{section}{Index}

\printindex

\end{document}